\documentclass[12pt]{amsart}
\usepackage{amssymb, amsmath, amsfonts, amsthm,hyperref,graphicx,xcolor}
\usepackage[all]{hypcap}
\usepackage[hmargin=1 in, vmargin = 1 in]{geometry}
\usepackage{bbm}
\usepackage{tikz}

\pgfdeclarelayer{bg}    
\pgfsetlayers{bg,main}  

\newcommand{\newword}[1]{\textbf{#1}}

\newenvironment{sbm}
    {\left[ \begin{smallmatrix}
    }
    { 
     \end{smallmatrix} \right]
    }

\newcommand{\floorAve}[2]{\lfloor \tfrac{#1+#2}{2} \rfloor}
\newcommand{\ceilAve}[2]{\lceil \tfrac{#1+#2}{2} \rceil}

\newcommand{\inmin}{\wedge}
\newcommand{\inmax}{\vee}

\newcommand{\Hull}{\mathrm{Hull}}

\newcommand{\Llog}{L \textrm{-} \! \log}

\newcommand{\proj}{\mathrm{proj}}

\newcommand{\Finite}{\mathrm{Finite}}
\newcommand{\Supp}{\mathrm{Support}}

\newcommand{\RR}{\mathbb{R}}

\newcommand{\ZZ}{\mathbb{Z}}

\newcommand{\cS}{\mathcal{S}}

\newcommand{\cW}{\mathcal{W}}

\newcommand{\tilh}{\widetilde{h}}
\newcommand{\tilg}{\widetilde{g}}

\newcommand{\One}{\mathbbm{1}}

\newcommand{\SkepExt}{\mathrm{SkepExt}}
\newcommand{\DB}{\partial^{\scalebox{0.5}{$\nwarrow$}}}
\newcommand{\HB}{\partial^{\scalebox{0.5}{$\leftarrow$}}}
\newcommand{\VB}{\partial^{\scalebox{0.5}{$\uparrow$}}}
\newcommand{\CB}{\partial^{\includegraphics[width=0.1 in]{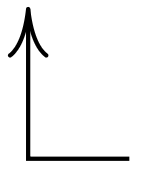}}}

\newcommand{\OB}{\partial^{+}}

\newcommand{\SHiveTop}{S_{\text{hive}}^{\text{top}}}
\newcommand{\SHiveBottom}{S_{\text{hive}}^{\text{bottom}}}
\newcommand{\SSkepTop}{S_{\text{skep}}^{\text{top}}}
\newcommand{\SSkepBottom}{S_{\text{skep}}^{\text{bottom}}}
\newcommand{\cSTop}{\cS^{\text{top}}}

\newcommand{\SWTri}[2]{\raisebox{-0.5\height}{${}_{#1 #2}$\!\! \includegraphics[width=0.25 in]{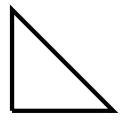}}}
\newcommand{\SETri}[2]{\raisebox{-0.5\height}{\includegraphics[width=0.25 in]{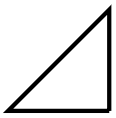}}_{#1 #2}}
\newcommand{\NWTri}[2]{\raisebox{0.1 in}{$\scriptstyle{#1#2}$}\raisebox{-0.5\height}{\includegraphics[width=0.25 in]{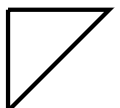}}} 
\newcommand{\NETri}[2]{\raisebox{-0.5\height}{\includegraphics[width=0.25 in]{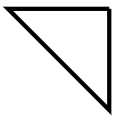}}^{#1 #2}}

\DeclareMathOperator{\zig}{\raisebox{-0.03 in}{\includegraphics[width=0.125 in]{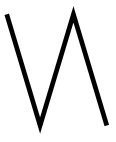}}}
\DeclareMathOperator{\zag}{\raisebox{-0.03 in}{\includegraphics[width=0.125 in]{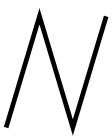}}}

\numberwithin{figure}{section}
\numberwithin{equation}{section}

\newtheorem{theorem}[equation]{Theorem}
\newtheorem{Theorem}[equation]{Theorem}
\newtheorem{lemma}[equation]{Lemma}
\newtheorem{Lemma}[equation]{Lemma}

\newtheorem{cor}[equation]{Corollary}
\newtheorem{Corollary}[equation]{Corollary}

\newtheorem*{LPP}{The LPP conjecture}
\newtheorem*{LPPMain}{Main Theorem}
\newtheorem*{corLPP}{Corollary~1}
\newtheorem*{AD}{The Ahlswede–Daykin inequality:}

\theoremstyle{definition}
\newtheorem{definition}[equation]{Definition}
\newtheorem{example}[equation]{Example}
\newtheorem{eg}[equation]{Example}
\newtheorem{remark}[equation]{Remark}

\newtheorem{Question}[equation]{Question}

\theoremstyle{plain}

\title[L-log-concavity and a proof of the LPP conjecture]{L-log-concavity and a proof of the conjecture of Lam, Postnikov and Pylyavskyy}
\author{David E Speyer}

\begin{document}
\begin{abstract}
Let $\lambda$, $\mu$, $\lambda'$, $\mu'$ be partitions. 
The conjecture of Lam, Postnikov and Pylyavskyy states that, if  $\lambda+\mu = \lambda' + \mu'$, and $\min(\lambda_i-\lambda_j, \mu_i-\mu_j) \leq \lambda'_i - \lambda'_j \leq \max(\lambda_i-\lambda_j, \mu_i-\mu_j)$ for all $1 \leq i<j \leq n$, then $s_{\lambda'} s_{\mu'} - s_{\lambda} s_{\mu}$ is Schur nonnegative. We prove this conjecture. 

Our proof is based on two key ideas. First, we introduce a new combinatorial model for Littlewood-Richardson coefficients which we name ``skeps", which are similar to but distinct from Knutson and Tao's hives. Second, we use tools from Murota's theory of $L$-convexity to prove an $\Llog$-concavity theorem for skeps.
\end{abstract}

\maketitle


A vector $\lambda = (\lambda_1, \lambda_2, \ldots, \lambda_n)$ in $\ZZ^n$ is called a \newword{partition} if $\lambda_1 \geq \lambda_2 \geq \cdots \geq \lambda_n \geq 0$.
For a partition $\lambda$, we write $s_{\lambda}$ for the Schur polynomial $s_{\lambda}$. A symmetric polynomial $f$ is called \newword{Schur nonnegative} if $f=\sum a_{\lambda} s_{\lambda}$ with $a_{\lambda} \geq 0$.

Given two vectors $x$ and $y \in \ZZ^n$, we define
\[ \Pi(x,y) := \{ z \in \ZZ^n : \min(x_i-x_j, y_i - y_j) \leq z_i - z_j \leq \max(x_i - x_j, y_i - y_j) \ \text{for all}\ 1 \leq i < j \leq n \} . \]
Note that, if $x'+y' = x+y$ and $x'$ is in $\Pi(x,y)$ then $y'$ is also in $\Pi(x,y)$.
See Lemma~\ref{PiGeometry} and Remark~\ref{PiGeometryRemark} for another description of $\Pi(x,y)$.

Lam, Postnikov and Pylyavskyy made the following conjecture, which we will prove:
\begin{LPP}
Let $\lambda$, $\mu$, $\lambda'$, $\mu'$ be partitions in $\ZZ^n$. Suppose that $\lambda'+\mu' = \lambda+\mu$ and $\lambda' \in \Pi(\lambda, \mu)$.  Then $s_{\lambda'}(x_1, x_2, \ldots, x_n) s_{\mu'}(x_1, x_2, \ldots, x_n) - s_{\lambda}(x_1, x_2, \ldots, x_n) s_{\mu}(x_1, x_2, \ldots, x_n)$ is Schur nonnegative. 
\end{LPP}

Lam, Postnikov and Pylyavskyy~\cite{LPP} proved many cases of this theorem but did not formally state this result as a conjecture; the LPP conjecture was formally stated as a conjecture by Dobrovolska and Pylyavskyy~\cite{DP}. 
Besides these papers, special cases of the LPP conjecture were proved in~\cite{KWvW, PvW, MvW, BallantineOrellana, FarberHopkinsTrongsiriwat, NguyenNguyenWoodruff}.
There is also significant literature proving that quantities related to Schur polynomials are $\log$-concave when their inputs are restricted to various arithmetic progressions, a condition which is much weaker than $\Llog$ concavity. (See Remark~\ref{Progression}.)
For example,~\cite{Oko1} proved monomial positivity of differences of the form $s_{k \lambda}^2 - s_{(k-1)\lambda} s_{(k+1)\lambda}$ and~\cite{HMMS} proved $\log$-concavity of Kostka numbers along arithmetic progressions in $(\ZZ^n)^2$. 

The \newword{Littlewood-Richardson coefficients} $c_{\lambda \mu}^{\nu}$ are defined by $s_{\lambda} s_{\mu} = \sum_{\nu} c_{\lambda \mu}^{\nu} s_{\nu}$. So the LPP conjecture says:
\begin{LPPMain} 
Let $\lambda$, $\mu$, $\lambda'$, $\mu'$ and $\nu$ be partitions in $\ZZ^n$. Suppose that $\lambda'+\mu' = \lambda+\mu$ and $\lambda' \in \Pi(\lambda, \mu)$.   Then $c_{\lambda' \mu'}^{\nu} \geq c_{\lambda \mu}^{\nu}$. 
\end{LPPMain}

We note that we have deliberately not said that, if  $\lambda'+\mu' = \lambda+\mu$ and $\lambda' \in \Pi(\lambda, \mu)$, then $s_{\lambda'} s_{\mu'} - s_{\lambda} s_{\mu}$ is Schur nonnegative, without writing the variables $(x_1, x_2, \ldots, x_n)$. 
That's because the notation $s_{\lambda}$, without explicit variables, denotes the Schur polynomial in an arbitrary large number of variables.  
The result is not true when interpreted that way: Take $n=1$, $(\lambda, \mu, \lambda', \mu') = ((1), (1), (2), (0))$. Then $\lambda+\mu=\lambda'+\mu'$ and $\lambda' \in \Pi(\lambda, \mu)$, but $s_2 s_0 = s_2 \not\succeq_s s_1^2 = s_2+s_{11}$. 
However, the LPP conjecture is still true, because $s_{11}(x_1) = 0$, so we have $s_2(x_1) s_0(x_1) = s_1(x_1)^2$. 

Here is a result to address the issue of adding additional variables.
\begin{corLPP} \label{PaddingCorrected}
Let $\lambda$, $\mu$, $\lambda'$, $\mu'$ and $\nu$ be partitions in $\ZZ^n$ with $\lambda+\mu = \lambda' + \mu'$,  $\lambda' \in \Pi(\lambda, \mu)$ and $\min(\lambda_i, \mu_i) \leq \lambda'_i \leq \max(\lambda_i, \mu_i)$ for $1 \leq i \leq n$. Then $s_{\lambda'} s_{\mu'} - s_{\lambda} s_{\mu}$ is Schur nonnegative.
\end{corLPP}

In Murota's language~\cite{Murota}, the conditions  $\lambda' \in \Pi(\lambda, \mu)$ and $\min(\lambda_i, \mu_i) \leq \lambda'_i \leq \max(\lambda_i, \mu_i)$ say that $\lambda'$ is in the $L^{\natural}$ convex hull of $\lambda$ and $\mu$. 
We will not use $L^{\natural}$-convexity elsewhere in this paper.

\begin{proof}
Given a vector $v = (v_1, v_2, \ldots, v_n)  \in \ZZ^n$, let $(v, 0^k):=(v_1, v_2, \ldots, v_n, 0,0,\ldots,0)$.
To say that $s_{\lambda'} s_{\mu'} - s_{\lambda} s_{\mu}$ is Schur nonnegative is to say that $c_{\lambda' \mu'}^{\nu} \geq c_{\lambda \mu}^{\nu}$ for $\nu \in \ZZ^{n+k}$, for any $k \geq 0$. 
Using the Main Theorem, this follows if $(\lambda', 0) \in \Pi((\lambda, 0^k), (\mu, 0^k))$. 
The required inequalities are thus 
\[ \begin{array}{ll}
\min(\lambda_i-\lambda_j, \mu_i-\mu_j) \leq \lambda'_i - \lambda'_j \leq \max(\lambda_i-\lambda_j, \mu_i-\mu_j) & 1 \leq i,j \leq n \\
 \min(\lambda_i-0, \mu_i-0) \leq \lambda'_i - 0 \leq \max(\lambda_i-0, \mu_i-0)  & 1\leq i \leq n .
 \end{array} \]
The first family of inequalities is the hypothesis $\lambda' \in \Pi(\lambda, \mu)$, the second family is the hypothesis $\min(\lambda_i, \mu_i) \leq \lambda'_i \leq \max(\lambda_i, \mu_i)$.
\end{proof}

\section{$L$-convexity and $\Llog$-concavity} \label{LConvex1}
Let $f$ be a function $\ZZ^n \to \RR \cup \{ \infty \}$. 
Murota~\cite{Murota} has developed an intricate theory of discrete convexity, describing ways in which $f$ is like a convex function on $\RR^n$.
Definitions~\ref{LConvexSetDefn} and~\ref{LConvexFunctionDefn} are equivalent to definitions of Murota;  we will check this equivalence in Section~\ref{LConvex2}.

\begin{definition}
Let $\One_n$ denote the vector $(1,1,\ldots, 1)$. 
\end{definition}

\begin{definition} \label{LConvexSetDefn}
Let $K$ be a subset of $\ZZ^n$. We will say that $K$ is \newword{$L$-convex} if, for any $x$, $y \in K$, we have $\Pi(x,y) \subseteq K$.
\end{definition}

 Note that $\Pi(x,y)$ is invariant under translation by $\One_n$, so all $L$-convex sets are invariant under translation by $\One_n$.
 We also remark explicitly that the empty set is $L$-convex.

\begin{definition} \label{LConvexFunctionDefn}
Let $f : \ZZ^n \to \RR \cup \{ \infty \}$ be a function. We say that $f$ is \newword{$L$-convex} if
\begin{enumerate}
\item If $x$, $y$, $x'$, $y' \in \ZZ^n$ with $x+y=x'+y'$ and $x' \in \Pi(x,y)$, then $f(x)+f(y) \geq f(x') + f(y')$.
\item There is a constant $r \in \RR$ such that $f(x + \One_n) = f(x)+r$ for all $x \in \ZZ^n$.
\end{enumerate}
\end{definition}

We make the following natural extensions of these definitions: 
\begin{definition}
A function $f : \ZZ^n \to \RR \cup \{ -\infty \}$ is \newword{$L$-concave} if $-f$ is $L$-convex.
\end{definition}

\begin{definition} \label{DefnLLogConcave}
A function $g: \ZZ^n \to \RR_{\geq 0}$ is \newword{$\Llog$-concave} if $\log g$ is $L$-concave. Here $\log 0$ is taken to be $-\infty$.
To be explicit, $g$ is $\Llog$-concave if and only if 
\begin{enumerate}
\item For any $x$, $y$, $x'$, $y' \in \ZZ^n$ with $x+y=x'+y'$ and $x' \in \Pi(x,y)$, we have $f(x') f(y') \geq f(x) f(y)$. 
\item There is a constant $R \in \RR_{>0}$ such that $f(x+\One_n) = R f(x)$ for all $x \in \ZZ^n$.
\end{enumerate}
\end{definition}

The first condition in Definition~\ref{DefnLLogConcave} is obviously reminiscent of the LPP conjecture. 
One might try to develop a general theory of Schur $\Llog$-concavity for functions valued in symmetric polynomials.
The author thinks this is a lovely idea, but he has not carried it out, we will discuss some of the issues in Remarks~\ref{SchurLlogConcaveRemark} and~\ref{2DoesNotImply1}.
Instead, this paper develops new ways of analyzing Littlewood-Richardson numbers and combines them with $\Llog$-concavity for ordinary $\RR_{\geq 0}$-valued functions, to prove the Main Theorem.
We now outline our proof.

\section{Overview of the proof}
We now describe a new combinatorial model for Littlewood-Richardson coefficients, which is at the heart of our proof.
Define
\[ 
\begin{array}{rcl}
\Delta_n &:=& \{ (i,j) \in \ZZ_{\geq 0}^2 : i+j \leq n \} \\
\Delta^+_n &:=& \{ (i,j) \in \Delta_n : i+j \equiv n \bmod 2\} \\
\Delta^-_n &:=& \{ (i,j) \in \Delta_n : i+j \equiv n-1 \bmod 2\} \\
\end{array} .\]

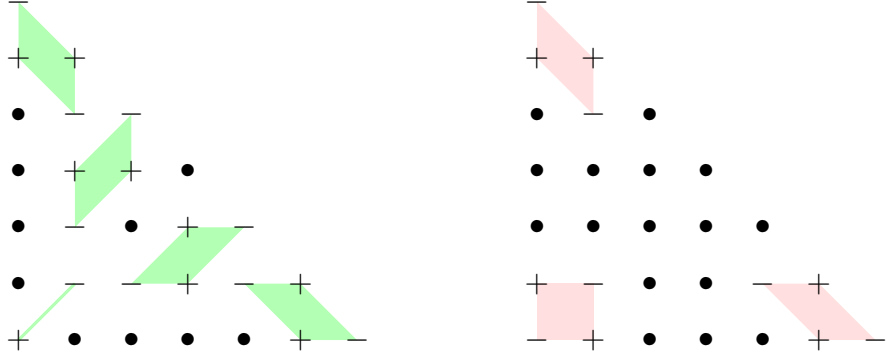
\begin{figure}
\begin{tikzpicture}[scale=0.75]
\node at (0,6) {$-$};
\node at (0,5) {$+$}; \node at (1,5) {$+$};
\node at (0,4) {$\bullet$}; \node at (1,4) {$-$}; \node at (2,4) {$-$};
\node at (0,3) {$\bullet$}; \node at (1,3) {$+$}; \node at (2,3) {$+$}; \node at (3,3) {$\bullet$};
\node at (0,2) {$\bullet$}; \node at (1,2) {$-$}; \node at (2,2) {$\bullet$}; \node at (3,2) {$+$}; \node at (4,2) {$-$};
\node at (0,1) {$\bullet$}; \node at (1,1) {$-$}; \node at (2,1) {$-$}; \node at (3,1) {$+$}; \node at (4,1) {$-$}; \node at (5,1) {$+$};
\node at (0,0) {$+$}; \node at (1,0) {$\bullet$}; \node at (2,0) {$\bullet$}; \node at (3,0) {$\bullet$}; \node at (4,0) {$\bullet$}; \node at (5,0) {$+$}; \node at (6,0) {$-$};

\begin{pgfonlayer}{bg}    
     \path [fill=green!30] (0,6) -- (1,5) -- (1,4) -- (0,5) -- (0,6);
      \path [fill=green!30] (2,4) -- (1,3) -- (1,2) -- (2,3) -- (2,4);
      \path [fill=green!30] (4,2) -- (3,1) -- (2,1) -- (3,2) -- (4,2);
     \path [fill=green!30] (6,0) -- (5,1) -- (4,1) -- (5,0) -- (6,0);
     \draw[green!30, very thick] (0,0) -- (1,1) ;
\end{pgfonlayer}
\end{tikzpicture} \qquad \qquad 
\begin{tikzpicture}[scale=0.75]
\node at (0,6) {$-$};
\node at (0,5) {$+$}; \node at (1,5) {$+$};
\node at (0,4) {$\bullet$}; \node at (1,4) {$-$}; \node at (2,4) {$\bullet$};
\node at (0,3) {$\bullet$}; \node at (1,3) {$\bullet$}; \node at (2,3) {$\bullet$}; \node at (3,3) {$\bullet$};
\node at (0,2) {$\bullet$}; \node at (1,2) {$\bullet$}; \node at (2,2) {$\bullet$}; \node at (3,2) {$\bullet$}; \node at (4,2) {$\bullet$};
\node at (0,1) {$+$}; \node at (1,1) {$-$}; \node at (2,1) {$\bullet$}; \node at (3,1) {$\bullet$}; \node at (4,1) {$-$}; \node at (5,1) {$+$};
\node at (0,0) {$-$}; \node at (1,0) {$+$}; \node at (2,0) {$\bullet$}; \node at (3,0) {$\bullet$}; \node at (4,0) {$\bullet$}; \node at (5,0) {$+$}; \node at (6,0) {$-$};

\begin{pgfonlayer}{bg}    
     \path [fill=pink!50] (0,6) -- (1,5) -- (1,4) -- (0,5) -- (0,6);
          \path [fill=pink!50] (0,0) -- (1,0) -- (1,1) -- (0,1) -- (0,0);
     \path [fill=pink!50] (6,0) -- (5,1) -- (4,1) -- (5,0) -- (6,0);
\end{pgfonlayer}
\end{tikzpicture}
\caption{The left hand side depicts the skep inequalities: In each green parallelogram, and in all translates thereof, the sum of the $+$ vertices is more than the sum of the $-$ vertices; we also impose this condition on the green line segment. The right hand side, with the pink parallelograms, depicts the hive inequalities in the same manner.} \label{SkepVersusHive}
\end{figure}

\begin{definition}
A \newword{skep} is a function $g : \Delta_n \to \ZZ$ obeying the following inequalities for all $(i,j)$ such that the indices remain within $\Delta_n$. 
\begin{equation}
\begin{array}{lcl}
g_{ij} - g_{(i+1)(j+1)} &\leq& g_{(i+1)j} - g_{(i+2)(j+1)} \\
g_{ij} - g_{(i+1)(j+1)} &\leq& g_{i(j+1)} - g_{(i+1)(j+2)} \\
g_{ij} - g_{(i+1)(j-1)} &\leq& g_{i(j-1)} - g_{(i+1)(j-2)} \\
g_{ij} - g_{(i+1)(j-1)} &\leq& g_{(i+1)j} - g_{(i+2)(j-1)} \\
0 &\leq& g_{00} - g_{11} \\
\end{array} \label{SkepInequalities}
\end{equation} 
\end{definition} 
See the left hand side of Figure~\ref{SkepVersusHive} for a visual representation of these inequalities. 

\begin{example} \label{IntroduceExample}
As our running example of a skep, we will take 
\[ 
g=
\begin{bmatrix}
11 &&&& \\
11 & 10 &&& \\
11 & 10 & 8 && \\
11 & 9 & 8 & 5 & \\
10 & 9 & 7 & 4 & 0 \\
\end{bmatrix} . \]
We are using Cartesian coordinates, so $g_{00}$ is the $10$ in the lower left, $g_{40}$ is the $0$ in the lower right and $g_{04}$ is the $11$ in the upper left. 
\end{example}

Skeps are related to, but not the same as, hives, which are functions $h: \Delta_n \to \ZZ$, obeying a different list of inequalities. 
(See Equation~\eqref{HiveInequalities} and Figure~\ref{SkepVersusHive}.)
The word ``skep" refers to an artificial beehive woven out of wicker, a tool used in medieval beekeeping. 
As we will see, our skeps are woven out of two parts.

\begin{remark} \label{NEinequalities}
The first two inequalities of~\eqref{SkepInequalities} say that the differences $g_{ij} - g_{(i+1)(j+1)}$ form an increasing tableaux, with $g_{00} - g_{11}$ smaller than all the others. 
Therefore, the inequality $g_{00} - g_{11} \geq 0$ forces all of the $g_{ij}-g_{(i+1)(j+1)}$ to be $\geq 0$.
\end{remark}

\begin{remark} \label{NWinequalities}
The third and fourth inequalities of~\eqref{SkepInequalities} say that the differences $g_{(i-1)(j+1)} - g_{ij}$ form a Gelfand-Tsetlin pattern, with smallest element $g_{0n} - g_{1(n-1)}$ and largest element $g_{(n-1)1} - g_{n0}$. 
\end{remark}

We now describe how skeps compute Littlewood-Richardson coefficients:

\begin{definition}
Given a function $g : \Delta_n \to \ZZ$, we define the \newword{diagonal boundary of $g$} to be
\[ \DB(g) := (g_{(n-1)1}-g_{n0}, g_{(n-2)2} - g_{(n-1)1},  \ldots, g_{0n} - g_{1(n-1)}) \]
and we define the \newword{corner boundary of $g$} to be 
\[ \CB(g) = (g_{(n-1)0} - g_{n0}, g_{(n-2)0} - g_{(n-1)0}, \ldots, g_{00} - g_{10}, g_{01} - g_{00}, \ldots, g_{0(n-1)} - g_{0(n-2)}, g_{0n}-g_{0(n-1)}) . \]
The shape of the arrow is meant to indicate the path that we travel through $\Delta_n$.
\end{definition}

Here is our big theorem saying that skeps compute Littlewood-Richardson coefficients. This theorem, like several in this section, is numbered according to where the reader can find its proof.

{
\renewcommand{\theequation}{\ref{SkepLR}}
\begin{Theorem} 
Let $\lambda$, $\mu$ and $\nu$ be partitions. Then $c_{\lambda \mu}^{\nu}$ is the number of skeps with $g_{n0}=0$, $\DB(g) = (\nu_1, \nu_2, \ldots, \nu_n)$ and $\CB(g) = (\lambda_1, \mu_1, \lambda_2, \mu_2, \ldots, \lambda_n, \mu_n)$.
\end{Theorem}
\addtocounter{equation}{-1}
}

\begin{example} \label{SkepBoundaryExample}
For the skep $g$ in Example~\ref{IntroduceExample}, we have
\begin{equation} \DB(g)=(5-0,\ 8-5,\ 10-8,\ 11-10) = (5, 3, 2, 1) \label{DBDefn} \end{equation}
and
\begin{equation} 
\begin{array}{rcl}
\CB(g) &=&  (4-0 ,\ 7-4,\ 9-7,\ 10-9,\ 11-10,\ 11-11,\ 11-11,\ 11-11) \\ &=& (4,3,2,1,1,0,0,0) .
\end{array}
 \label{CBDefn} \end{equation}
This is one of the skeps that computes the Littlewood-Richardson coefficient $c_{4210,\ 3100}^{5321}$. 
\end{example}


\begin{definition}
Given a function $g : \Delta_n \to \ZZ$, we define $g^+$ and $g^-$ to be the restrictions $g|_{\Delta_n^+}$ and $g|_{\Delta_n^-}$. 
Conversely, given $g^+ : \Delta_n^+ \to \ZZ$ and $g^- : \Delta_n^- \to \ZZ$, we define $(g^+, g^-) : \Delta_n \to \ZZ$ to be the function with $(g^+, g^-)|_{\Delta_n^+} = g^+$ and $(g^+, g^-)|_{\Delta_n^-} = g^-$.
\end{definition}

\begin{example} \label{InterwovenExample}
For the skep in Example~\ref{IntroduceExample}, we have
\[ 
g^+=
\begin{bmatrix}
11 &&&& \\
 & 10 &&& \\
11 &  & 8 && \\
 & 9 &  & 5 & \\
10 &  & 7 & & 0 \\
\end{bmatrix} 
 \ \text{and}\ g^- = 
 \begin{bmatrix}
 &&&& \\
11 &  &&& \\
 & 10 &  && \\
11 &  & 8 &  & \\
& 9 &  & 4 &  \\
\end{bmatrix} 
 . \]
\end{example}

Note that, in each of the inequalities in~\eqref{SkepInequalities}, one side involves only terms from $g^+$ and the other side involves only terms from $g^-$. 
We think of $g^+$ and $g^-$ as the two halves of the skep, and we think of~\eqref{SkepInequalities} as describing when we can interweave them. 

\begin{definition} \label{BoundaryParts}
Given $g : \Delta^n \to \ZZ$, let $\CB(g) = (z_1, z_2, \ldots, z_{2n-1}, z_{2n})$. We define
\[
 \begin{array}{lcl}
\partial_1(g) &=& (z_1, z_3, z_5, \ldots, z_{2n-1}) \\
\partial_2(g) &=& (z_2, z_4, z_6, \ldots, z_{2n}) \\
\partial^+(g) &=& \partial_1(g)+\partial_2(g) = (z_1+z_2, z_3+z_4, \ldots, z_{2n-1}+z_{2n}) \\
\HB(g) &=& (z_1, z_2, \ldots, z_n) \\
\VB(g) &=& (z_{n+1}, z_{n+2}, \ldots, z_{2n}) \\
\end{array}\] 
\end{definition}

In this language, Theorem~\ref{SkepLR} says that $c_{\lambda \mu}^{\nu}$ is the number of skeps with $(\partial_1(g), \partial_2(g), \DB(g)) = (\lambda, \mu, \nu)$ and $g_{n0} = 0$.

\begin{definition}
Note that Equation~\eqref{DBDefn}, defining  $\DB(g)$, only uses terms from $g^+$. Thus, for any $g^+ : \Delta_n^+ \to \ZZ$, we use Equation~\eqref{DBDefn} to define $\DB(g^+)$. 

Likewise, note that
\begin{equation}
\OB(g) = (g_{(n-2)2} - g_{n0}, g_{(n-4)4} - g_{(n-2)2} , \ldots, g_{0n} - g_{2(n-2)}) \label{PlusDefn} 
 \end{equation}
 only uses terms from $g^+$. Thus, for any $g^+ : \Delta_n^+ \to \ZZ$, we use Equation~\eqref{PlusDefn} to define $\OB(g^+)$. 
 We call $\OB(g^+)$ the \newword{outer boundary} of $g^+$.
\end{definition}

\begin{definition}
Given $g^+ : \Delta^+_n \to \ZZ$ and a dominant vector $\lambda$, we define $\SkepExt(g^+, \lambda)$ to be the number of $g^-: \Delta^-_n \to \ZZ$ so that $(g^+, g^-)$ is a skep with $\partial_1(g^+, g^-) = \lambda$. 
\end{definition}

We can thus rewrite Theorem~\ref{SkepLR} as follows:
{
\renewcommand{\theequation}{\ref{SkepLRAsSum}}
\begin{Corollary} 
Let $\lambda$, $\mu$ and $\nu$ be partitions. Set $\pi = \lambda+\mu$.
Then
\[
c_{\lambda \mu}^{\nu} = \sum \SkepExt(g^+, \lambda)
\]
where the sum runs over $g^+ : \Delta_n^+ \to \ZZ$ such that $\DB(g^+) = \nu$ and  $\OB(g^+) = \pi$.
\end{Corollary}
\addtocounter{equation}{-1}
}

\begin{remark}
We described the sum in Corollary~\ref{SkepLRAsSum} as running over an infinite set of possible $g^+$. 
We can restrict the sum to a finite set: If $(g^+, g^-)$ is a skep, then Remark~\ref{NEinequalities} shows that
\[ g^+_{k(n-k)} \leq g^+_{(k-1)(n-k-1)} \leq g^+_{(k-2)(n-k-2)} \leq \cdots \leq 
\begin{cases} g_{0(n-2k)} & k \leq n/2 \\ g_{(2k-n)0} & k \geq n/2 \end{cases} . \]
Since $g^+_{k(n-k)}$,  $g^+_{i0}$ and $g^+_{0j}$ are fixed by the conditions that $\DB(g^+) = \nu$ and  $\OB(g^+) = \pi$, this 
restricts the remaining $g^+_{ij}$ values to lie in a bounded range. 
\end{remark}
%

Thus, the Main Theorem will follow from:
{
\renewcommand{\theequation}{\ref{BetterLPP}}
\begin{Theorem} 
Let $\lambda$, $\mu$ and $\nu$ be partitions. Set $\pi = \lambda+\mu$.
Fix $g^+ : \Delta^+_n \to \ZZ$ with $\DB(g^+) = \nu$ and $\OB(g^+) = \pi$.
Let $\lambda'+\mu' = \pi$ with $\lambda' \in \Pi(\lambda, \mu)$. 
Then
\[
\SkepExt(g^+, \lambda)  \leq \SkepExt(g^+, \lambda') .
\]
\end{Theorem}
\addtocounter{equation}{-1}
}

The main (but not only) tool in proving Theorem~\ref{BetterLPP} is the following result, which explains why we have spent so much time on $\Llog$-concavity.
{
\renewcommand{\theequation}{\ref{SkepExtLConcave}}
\begin{Theorem}
Fix $g^+ : \Delta^+_n \to \ZZ$. Then $\SkepExt(g^+, \lambda)$ is $\Llog$-concave as a function of $\lambda$.
\end{Theorem}
\addtocounter{equation}{-1}
}

A sum of $\Llog$-concave functions need not be $\Llog$-concave. 
Therefore, Theorem~\ref{SkepExtLConcave} does \textbf{not} imply that the function $\lambda \mapsto c_{\lambda (\pi-\lambda)}^{\nu}$ is $\Llog$-concave. 
We therefore pose the question:
\begin{Question}
Fix $\nu$ and $\pi$ with $|\nu| = |\pi|$. Is the function $\lambda \mapsto c_{\lambda (\pi-\lambda)}^{\nu}$  $\Llog$-concave?
\end{Question}
%

The proof of Theorem~\ref{SkepExtLConcave} uses very little of the structure of skeps. 
Instead, our proof is based on the following result, which refers only to general concepts of $\Llog$-concavity.

{
\renewcommand{\theequation}{\ref{PointCountLConcave}}
\begin{theorem} 
Let $K \subset \ZZ^{M+N}$ be an $L$-convex set. For $x \in \ZZ^M$, let \[E_K(x) = \#\{ y \in \ZZ^N : (x,y) \in K \},\]and assume that $E_K(x) < \infty$ for all $x \in \ZZ^M$. Then $E_K$ is an $\Llog$-concave function on $\ZZ^M$.
\end{theorem}
\addtocounter{equation}{-1}
}

This is a good point to give an overview of the remainder of the paper.
Section~\ref{HiveSection} proves Theorem~\ref{SkepLR} and Corollary~\ref{SkepLRAsSum}.
Section~\ref{LConvex2} gives more background on $L$-convexity.
Section~\ref{LConcaveProof}  proves Theorem~\ref{PointCountLConcave}.
In Section~\ref{Finishing}, we apply these results to prove  Theorems~\ref{SkepExtLConcave}, \ref{BetterLPP} and the Main Theorem.
In Section~\ref{MinimalInequalities}, we give a list of $(\lambda, \mu, \lambda', \mu')$ such that proving the LPP conjecture for those quadruples establishes the full conjecture.


\section{Hives, the octahedron recurrence, and proof of Theorem~\ref{SkepLR}} \label{HiveSection}

Skeps are based on hives, which were invented by Knutson and Tao~\cite[Appendix 1]{KT1}. See also~\cite{Buch}.
Like a skep, a hive is a triangular array of integers, obeying certain inequalities, and are used to compute Littlewood-Richardson coefficients. 
However, the details are different.

\begin{definition}
A \newword{hive} is a function $h : \Delta_n^+ \to \ZZ$  obeying the following inequalities for all $(i,j)$ such that the indices remain within $\Delta_n$. See the right hand side of Figure~\ref{SkepVersusHive}.
\begin{equation}
 \begin{array}{lcl}
h_{(i+1)j} + h_{i(j+1)} &\geq& h_{ij} + h_{(i+1)(j+1)} \\
h_{i(j+1)} + h_{ij} &\geq& h_{(i-1)(j+1)} + h_{(i+1)j} \\
h_{(i+1)j} + h_{ij} &\geq& h_{(i+1)(j-1)} + h_{i(j+1)} \\
\end{array} \label{HiveInequalities}
\end{equation} 
\end{definition}

\begin{example} \label{HiveExample}
Our running example of a hive will be
\[ h = \begin{bmatrix}
11 &&&& \\ 
11&10 &&&\\
11&10&8 && \\
10&10&8&5& \\
7&7&6&4&0 \\
\end{bmatrix} . \]
\end{example}


Here is the Theorem of Knutson and Tao~\cite[Appendix 1]{KT1}. See also~\cite{KTW} for a direct proof.
\begin{Theorem}  \label{HiveLR}
Let $\lambda$, $\mu$ and $\nu$ be partitions. Then $c_{\lambda \mu}^{\nu}$ is the number of hives with $h_{n0}=0$, $\HB(h) = \lambda$, $\VB(h) = \mu$ and $\DB(h) = \nu$. (See Definition~\ref{BoundaryParts} for these notations.)
\end{Theorem}

\begin{example}
For the hive $h$ in Example~\ref{HiveExample}, we have
\[ 
\begin{array}{lclcl}
\HB(h) &=& (4-0, 6-4, 7-6, 7-7) &=& (4,2,1,0) \\
\VB(h) &=& (10-7, 11-10, 11-11, 11-11) &=& (3,1,0,0) \\
\DB(h)  &=& (5-0, 8-5, 10-8, 11-10) &=& (5,3,2,1) .
  \end{array}\]
\end{example}

Littlewood-Richardson coefficients are symmetric in $\lambda$ and $\mu$, but the definition of hives is not symmetric.
Henriques and Kamnitzer~\cite{HK} found a subtle piecewise linear bijection between the hives with boundary $(\lambda, \mu, \nu)$  and the hives with boundary $(\mu, \lambda, \nu)$, using the \newword{octahedron recurrence}.
Our presentation largely follows \cite[Section 3]{HK}; see also~\cite{KTW}.

\begin{definition}
Define 
\[ T_n = \{ (i,j,t) : 0 \leq i,j,\ |t| \leq n-i-j,\ t \equiv i+j+n \bmod 2 \}. \]
 $T_n$ is the intersection of a tetrahedron with a face centered cubic lattice.
\end{definition}

\begin{definition}
A function $\tilh: T_n \to \ZZ$ is said to obey the \newword{octahedron recurrence} if, for all $(i,j,t) \in T_n$ such that $(i,j,t-2)$ also lies in $T_n$, we have
\begin{equation} \tilh(i,j,t) + \tilh(i,j,t-2) \! = \!
\left\{ \! \! \begin{array}{l@{\ }l} 
\max\left( \begin{matrix} \tilh(i-1,j,t-1) + \tilh(i+1,j,t-1), \quad \\ \quad \tilh(i,j-1,t-1) + \tilh(i,j+1,t-1) \end{matrix} \right) & i,j \geq 1 \\[0.5 cm]
\tilh(i-1,0,t-1) + \tilh(i+1,0,t-1) & i \geq 1,\ j=0 \\[0.3 cm]
 \tilh(0,j-1,t-1) + \tilh(0,j+1,t-1) & i=0,\ j \geq 1 \\[0.3 cm]
 \tilh(1,0,t-1)+\tilh(0,1,t-1) & i=j=0 \\
\end{array} \right.
\label{OctRecurr}
 \end{equation}
\end{definition}

\begin{example} \label{OctExample}
We will present an example of a function obeying the octahedron recurrence, for $n=4$, by drawing slices at each value of $t$. 
To visualize them in three dimensions, stack them up, aligning the lower left corners. 
\[ \begin{array}{cccccccccc}
\raisebox{\height}{\text{slice:}} &
\raisebox{\height}{$\begin{sbm}   7 \end{sbm}$} &
\raisebox{\height}{$\begin{sbm}   10& \\ &7 \end{sbm}$} &
\raisebox{\height}{$\begin{sbm}   11&& \\ &10& \\ 10&&6 \end{sbm}$} &
\raisebox{\height}{$\begin{sbm}   11&&& \\ &10&& \\ 11&&8& \\ &9&&4 \end{sbm}$} &
\raisebox{\height}{$\begin{sbm} 11&&&& \\ &10&&& \\ 11&&8&& \\ &9&&5& \\ 10&&7&&0 \end{sbm}$} &
\raisebox{\height}{$\begin{sbm}   11&&& \\ &9&& \\ 10&&7& \\ &8&&3 \end{sbm}$} &
\raisebox{\height}{$\begin{sbm}   10&& \\ &8& \\ 8&&4 \end{sbm}$} &
\raisebox{\height}{$\begin{sbm}   8& \\ &4 \end{sbm}$} &
\raisebox{\height}{$\begin{sbm}   4 \end{sbm}$}  \\
t= & -4 & -3 & -2 & -1 & 0 & 1 & 2 & 3 & 4 \\ 
\end{array} \]
For example, the $7$ in the $t=1$ slice is $\max(9+5, 8+7) - 8$, where the values inside the $\max$ come from the $t=0$ slice and the $8$ which is subtracted off comes from the $t=-1$ slice.
\end{example}

Here is how Henriques and Kamnitzer use the octahedron recurrence to biject hives with boundary $(\lambda, \mu, \nu)$ and hives with boundary $(\mu, \lambda, \nu)$.
\begin{theorem} \label{MainHK}
Let $\tilh$ be a function on $T_n$ obeying the octahedron recurrence. Define functions $h : \Delta_n \to \ZZ$ and $h' : \Delta_n \to \ZZ$ by $h_{ij} := \tilh_{ij(i+j-n)}$ and $h'_{ij} := \tilh_{ij(n-i-j)}$. Then $h$ is a hive if and only if $h'$ is a hive. Moreover,  $\HB(h) = \VB(h')$, $\VB(h) = \HB(h')$ and $\DB(h) = \DB(h')$
\end{theorem}

Each of $h$ and $h'$ determine each other via the octahedron recurrence, so this is a bijection between hives with boundary $(\lambda, \mu, \nu)$ and hives with boundary $(\mu, \lambda, \nu)$. 

\begin{example}
Example~\ref{OctExample} was produced by this method. The hive $h$ from Example~\ref{HiveExample} was written on the diagonals of the slices at heights $-5 \leq t \leq 0$.  Then we can read $h'$ off the diagonals in the slices at heights $0 \leq t \leq 5$, giving 
\[ h' = \begin{bmatrix}
11&&&& \\
11&10&&& \\
10&9&8&& \\
8&8&7&5 \\
4&4&4&3&0 \\
\end{bmatrix} .\]
Note that this hive has boundary $(3100, 4210, 5321)$. 
\end{example}

The reader may recognize the slices at $t=0$ and $t=-1$ in Example~\ref{OctExample} -- they are $g^+$ and $g^-$ from Example~\ref{InterwovenExample}. 
This is the general story: If we start with a hive on the bottom of $T_n$ and run the octahedron recurrence up to the halfway point, we get a skep.
We now give the details.

\begin{definition} \label{HiveRestriction}
We write $\proj$ for the projection $\proj(i,j,t) = (i,j)$ from $T_n$ to $\Delta_n$. 
Let $\tilh$ be any function $T_n \to \ZZ$ and let $S$ be a subset of $T_n$ such that $\proj : S \to \Delta_n$ is bijective. 
We abuse notation slightly and define $\tilh|_S$ to be the composition $\Delta_n \xrightarrow{\proj^{-1}} S  \xrightarrow{\ \ \ \tilh \ \ \ }\ZZ$.
Thus, we say that  \newword{$\tilh|_S$ is a hive} or \newword{is a skep}  if this composition is a hive or a skep, and we use notation like $\HB(\tilh|_S)$ etcetera, which we have previously defined for functions $\Delta_n \to \ZZ$.
\end{definition}

We will want to apply Definition~\ref{HiveRestriction} to four particular sets $S$:
\begin{definition} \label{KeySlices}
Let $\epsilon(i,j)$ be $0$ if $(i,j) \in \Delta_n^+$ and $1$ if $(i,j) \in \Delta_n^-$.
We define
\[ \begin{array}{lcl}
\SHiveTop &=& \{ (i,j, n-i-j) : (i,j) \in \Delta_n \} \\[0.25 cm]
\SHiveBottom &=&   \{ (i,j, -n+i+j) : (i,j) \in \Delta_n \}\\[0.25 cm]
\SSkepTop &=& \{ (i,j,\epsilon(i,j)): (i,j) \in \Delta_n  \} \\[0.25 cm]
\SSkepBottom &=& \{ (i,j,-\epsilon(i,j)): (i,j) \in \Delta_n  \} . \\
\end{array} \]
Clearly, for each of these subsets of $T_n$, the projection down to $\Delta_n$ is bijective. 
\end{definition}

We now introduce vocabulary to describe various inequalities which we will impose on $\tilh$.

\begin{definition} \label{DefnTrianglesRhombiInequalities}
We define a \newword{unit triangle} to be an equilateral triangle whose vertices are in $T_n$ and whose edges have length $\sqrt{2}$. 
In other words, the vertices of a unit triangle are of the form $w\pm(1,0,0)$, $w \pm (0,1,0)$, $w \pm (0,0,1)$ for some lattice point $w$ and some (independent) choices of signs, such that all three vertices are in $T_n$.

We define a \newword{unit rhombus} to be two coplanar adjacent unit triangles, bordering along an edge. 
A unit rhombus has two $60^{\circ}$ angles and two $120^{\circ}$ angles.
We call the line segment between the $60^{\circ}$ angles the \newword{long diagonal}, and we call the line segment between the $120^{\circ}$ angles the \newword{short diagonal}.

Let $\tilh$ be a function $T_n \to \ZZ$ and let $R$ be a unit rhombus with long diagonal $(x,y)$, and short diagonal $(x', y')$. We say that $\tilh$ obeys the \newword{rhombus inequality at $R$} if $\tilh(x)+\tilh(y) \leq \tilh(x') + \tilh(y')$. We say that $\tilh$ \newword{obeys all rhombus inequalities} if $\tilh$ obeys the rhombus inequality at every unit rhombus in $T_n$.
\end{definition}


The main goal of this section is to show:
\begin{Theorem} \label{TechnicalOctTheorem}
Let $\tilh : T_n \to \ZZ$ be a function obeying the octahedron recurrence. Then the following are equivalent
\begin{enumerate}
\item $\tilh|_{\SHiveTop}$ is a hive.
\item  $\tilh|_{\SHiveBottom}$ is a hive.
\item $\tilh|_{\SSkepTop}$ is a skep.
\item $\tilh|_{\SSkepBottom}$ is a skep.
\item $\tilh$ obeys all rhombus inequalities.
\end{enumerate}
Moreover, their boundaries are related by:
\[ \DB(\tilh|_{\SHiveTop}) = \DB(\tilh|_{\SSkepTop}) = \DB(\tilh|_{\SSkepBottom}) = \DB(\tilh|_{\SHiveBottom}) \]
\[ \HB(\tilh|_{\SHiveTop}) = \partial_1(\tilh|_{\SSkepTop}) = \partial_2(\tilh|_{\SSkepBottom}) = \VB(\tilh|_{\SHiveBottom})  \]
\[ \VB(\tilh|_{\SHiveTop}) = \partial_2(\tilh|_{\SSkepTop}) = \partial_1(\tilh|_{\SSkepBottom}) = \HB(\tilh|_{\SHiveBottom}) . \]
\end{Theorem}
The equivalence of $(1)$, $(2)$ and $(5)$ is due to~\cite{HK}, as are the results relating the boundaries of $\tilh|_{\SHiveTop}$ and  $\tilh|_{\SHiveBottom}$. As we will see, $(3)$ and $(4)$ is also implicit in the work of~\cite{HK}, but require significant unpacking. 

We first deduce consequences of Theorem~\ref{TechnicalOctTheorem}, and then we prove it.

\begin{Theorem} \label{SkepLR}
Let $\lambda$, $\mu$ and $\nu$ be partitions. Then $c_{\lambda \mu}^{\nu}$ is the number of skeps with $g_{n0}=0$, $\DB(g) = (\nu_1, \nu_2, \ldots, \nu_n)$ and $\partial_1(g) = \lambda$ and $\partial_2(g) = \mu$.
\end{Theorem}

\begin{proof}[Proof of Theorem~\ref{SkepLR}, assuming Theorem~\ref{TechnicalOctTheorem}]
Theorem~\ref{TechnicalOctTheorem} gives a bijection between 
\begin{enumerate}
\item[(A)] hives with $h_{n0} = 0$, $\DB(h) = \nu$, $\HB(h) = \lambda$ and $\VB(h) = \mu$ and
\item[(B)] skeps with $g_{n0}=0$, $\DB(g) = \nu$, $\partial_1(g) = \lambda$ and $\partial_2(g) = \mu$. 
\end{enumerate} 
Specifically, given a hive as in~(A), place it on $\SHiveBottom$, use the octahedron recurrence to compute a function $\tilh$ on $T_n$, and consider the restriction $\tilh|_{\SSkepBottom}$. 
The inverse map, given a skep as in~(B), is to place it on $\SSkepBottom$, use the octahedron recurrence to compute a function $\tilh$ on $T_n$, and  consider the restriction $\tilh|_{\SHiveBottom}$. 

Theorem~\ref{HiveLR} tells us that the number of hives as in condition~(A) is $c_{\lambda \mu}^{\nu}$, so the same is true of the number of skeps obeying condition~(B).
\end{proof}

\begin{Corollary} \label{SkepLRAsSum}
Let $\lambda$, $\mu$ and $\nu$ be partitions. Set $\pi = \lambda+\mu$.
Then
\begin{equation}
c_{\lambda \mu}^{\nu} = \sum \SkepExt(g^+, \lambda) \label{SkepSum}
\end{equation}
where the sum runs over $g^+ : \Delta_n^+ \to \ZZ$ such that $\DB(g^+) = \nu$ and  $\OB(g^+) = \pi$.
\end{Corollary}

\begin{proof}
$\SkepExt(g^+,  \lambda)$ is the number of skeps $(g^+, g^-)$ with $\partial_1(g^+, g^-) = \lambda$.
Thus, the full list of conditions is that $(g^+, g^-)$ is a skep, $\OB(g^+) = \pi$, $\DB(g^+) = \nu$ and $\partial_1(g^+, g^-) = \lambda$.
Since $\OB(g^+) = \partial_1(g^+, g^-) + \partial_2(g^+, g^-)$ and $\pi = \lambda+\mu$, we may replace the condition  $\OB(g^+) = \pi$ by $\partial_2(g^+, g^-) = \mu$.
This is precisely the list of conditions in  Theorem~\ref{SkepLR}.
\end{proof}

Since $c_{\lambda \mu}^{\nu} = c_{\mu \lambda}^{\nu}$, Equation~\ref{SkepSum} tells us that
\[  \sum \SkepExt(g^+, \lambda)  =  \sum \SkepExt(g^+, \mu) \]
where the sum runs over $g^+ : \Delta_n^+ \to \ZZ$ such that $\DB(g^+) = \nu$ and  $\OB(g^+) = \lambda+\mu$.
We will now prove a refined version of that statement.

\begin{Theorem} \label{SkepCommutative}
Let $\lambda$ and $\mu$ be partitions. Let $g^+ : \Delta^+ \to \ZZ$ with $\OB(g^+) = \lambda+\mu$. Then
\[ \SkepExt(g^+, \lambda) = \SkepExt(g^+, \mu) .\]
\end{Theorem}

\begin{proof}[Proof of Theorem~\ref{SkepCommutative}, assuming Theorem~\ref{TechnicalOctTheorem}]
We must give a bijection between
\begin{enumerate}
\item[(C)] functions $g_1^- : \Delta_n^- \to \ZZ$ such that $(g^+, g_1^-)$ is a skep and $\partial_1(g^+, g_1^-) = \lambda$ and
\item[(D)] functions $g_2^- : \Delta_n^- \to \ZZ$ such that $(g^+, g_2^-)$ is a skep and $\partial_1(g^+, g_2^-) = \mu$.
\end{enumerate}
Let $g_1^-$ be as in~(C) and put $g = (g^+, g_1^-)$. Since $\partial_1(g) = \lambda$, we must have $\partial_2(g) = \OB(g^+) - \partial_1(g) = (\lambda+\mu) - \lambda = \mu$. 
Now, place the skep $g$ on $\SSkepBottom$, use the octahedron recurrence to propagate to a function $\tilg$ on $T_n$, and let $g'$ be the restriction of $\tilg$ to $\SSkepTop$. 
Since $\Delta^+_n \times \{ 0 \}$ is in both $\SSkepBottom$ and $\SSkepTop$, we have $(g')^+ = g^+$, so $g'$ is of the form $(g^+, g_2^-)$ for some $g_2^- : \Delta_n^- \to \ZZ$. 
By Theorem~\ref{TechnicalOctTheorem}, $g'$ is a skep, and $\partial_1(g') = \partial_2(g) = \mu$. So we have constructed $g_2^-$ as in (D). The inverse map, from~(D) to~(C), is analogous.
\end{proof}

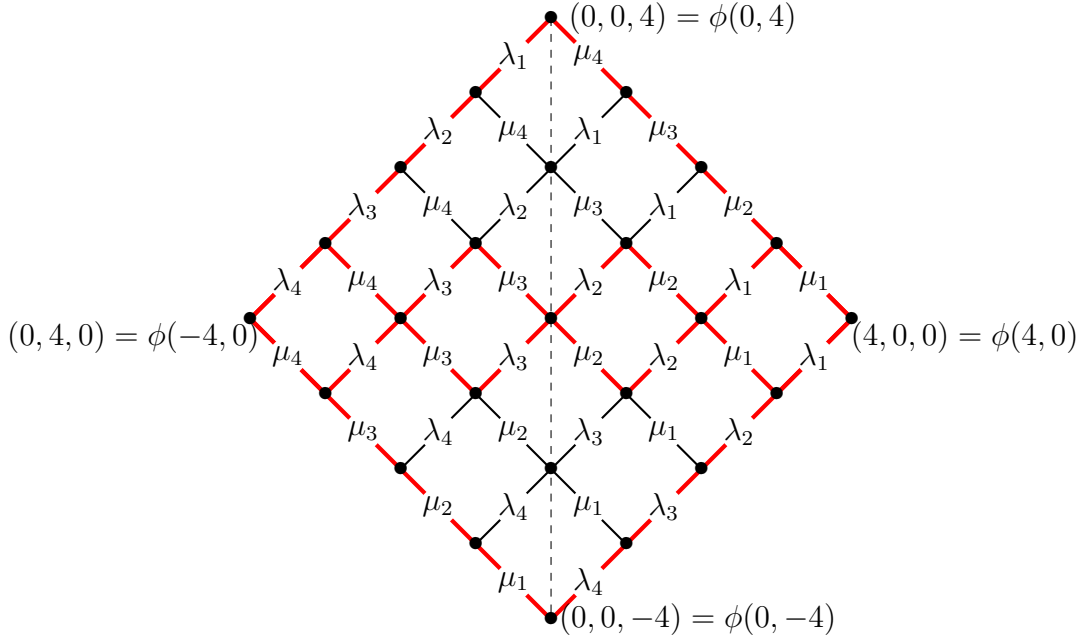
\begin{figure}
\centerline{\begin{tikzpicture}
\begin{pgfonlayer}{bg}    
\draw[thick] (0,4) -- (4,0);
\draw[thick] (-1,3) -- (3,-1);
\draw[thick] (-2,2) -- (2,-2);
\draw[thick] (-3,1) -- (1,-3);
\draw[thick] (-4,0) -- (0,-4);
\draw[thick] (0,4) -- (-4,0);
\draw[thick] (1,3) -- (-3,-1);
\draw[thick] (2,2) -- (-2,-2);
\draw[thick] (3,1) -- (-1,-3);
\draw[thick] (4,0) -- (0,-4);

\draw[ultra thick, red] (4,0) -- (0,-4) -- (-4,0);
\draw[ultra thick, red] (4,0) -- (0,4) -- (-4,0);
\draw[ultra thick, red] (4,0) -- (3,-1) -- (2,0) -- (1,-1) -- (0,0) -- (-1,-1) -- (-2, 0) -- (-3, -1) -- (-4, 0) ;
\draw[ultra thick, red] (4,0) -- (3,1) -- (2,0) -- (1,1) -- (0,0) -- (-1,1) -- (-2, 0) -- (-3, 1) -- (-4, 0) ; 

\draw[thin, dashed] (0,4) -- (0,-4);
\end{pgfonlayer}

\node at (0,4) {$\bullet$};
\node at (-1,3) {$\bullet$};\node at (1,3) {$\bullet$};
\node at (-2,2) {$\bullet$};\node at (0,2) {$\bullet$};\node at (2,2) {$\bullet$};
\node at (-3,1) {$\bullet$};\node at (-1,1) {$\bullet$};\node at (1,1) {$\bullet$};\node at (3,1) {$\bullet$};
\node at (-4,0) {$\bullet$};\node at (-2,0) {$\bullet$};\node at (0,0) {$\bullet$};\node at (2,0) {$\bullet$};\node at (4,0) {$\bullet$};
\node at (-3, -1) {$\bullet$};\node at (-1, -1) {$\bullet$};\node at (1, -1) {$\bullet$};\node at (3, -1) {$\bullet$};
\node at (-2, -2) {$\bullet$};\node at (0, -2) {$\bullet$};\node at (2, -2) {$\bullet$};
\node at (-1, -3) {$\bullet$};\node at (1, -3) {$\bullet$};
\node at (0, -4) {$\bullet$};

\node at (1.75, 4) {$(0,0,4) = \phi(0,4)$};
\node at (5.50, -0.25) {$(4,0,0) = \phi(4,0) $};
\node at (-5.55, -0.25) {$(0,4,0) = \phi(-4,0)$};
\node at (1.95, -4) {$(0,0,-4) = \phi(0,-4)$};

\fill[white] (3.5, -0.5) circle (0.25 cm);
\node at (3.5, -0.5) {$\lambda_1$};
\fill[white] (2.5, 0.5) circle (0.25 cm);
\node at (2.5, 0.5) {$\lambda_1$};
\fill[white] (1.5, 1.5) circle (0.25 cm);
\node at (1.5, 1.5) {$\lambda_1$};
\fill[white] (0.5, 2.5) circle (0.25 cm);
\node at (0.5, 2.5) {$\lambda_1$};
\fill[white] (-0.5, 3.5) circle (0.25 cm);
\node at (-0.5, 3.5) {$\lambda_1$};

\fill[white] (2.5, -1.5) circle (0.25 cm);
\node at (2.5, -1.5) {$\lambda_2$};
\fill[white] (1.5, -0.5) circle (0.25 cm);
\node at (1.5, -0.5) {$\lambda_2$};
\fill[white] (0.5, 0.5) circle (0.25 cm);
\node at (0.5, 0.5) {$\lambda_2$};
\fill[white] (-0.5, 1.5) circle (0.25 cm);
\node at (-0.5, 1.5) {$\lambda_2$};
\fill[white] (-1.5, 2.5) circle (0.25 cm);
\node at (-1.5, 2.5) {$\lambda_2$};

\fill[white] (1.5, -2.5) circle (0.25 cm);
\node at (1.5, -2.5) {$\lambda_3$};
\fill[white] (0.5, -1.5) circle (0.25 cm);
\node at (0.5, -1.5) {$\lambda_3$};
\fill[white] (-0.5, -0.5) circle (0.25 cm);
\node at (-0.5, -0.5) {$\lambda_3$};
\fill[white] (-1.5, 0.5) circle (0.25 cm);
\node at (-1.5, 0.5) {$\lambda_3$};
\fill[white] (-2.5, 1.5) circle (0.25 cm);
\node at (-2.5, 1.5) {$\lambda_3$};

\fill[white] (0.5, -3.5) circle (0.25 cm);
\node at (0.5, -3.5) {$\lambda_4$};
\fill[white] (-0.5, -2.5) circle (0.25 cm);
\node at (0-.5, -2.5) {$\lambda_4$};
\fill[white] (-1.5, -1.5) circle (0.25 cm);
\node at (-1.5, -1.5) {$\lambda_4$};
\fill[white] (-2.5, -0.5) circle (0.25 cm);
\node at (-2.5, -0.5) {$\lambda_4$};
\fill[white] (-3.5, 0.5) circle (0.25 cm);
\node at (-3.5, 0.5) {$\lambda_4$};

\fill[white] (3.5,0.5) circle (0.25 cm);
\node at (3.5, 0.5) {$\mu_1$};
\fill[white] (2.5, -0.5) circle (0.25 cm);
\node at (2.5, -0.5) {$\mu_1$};
\fill[white] (1.5, -1.5) circle (0.25 cm);
\node at (1.5, -1.5) {$\mu_1$};
\fill[white] (0.5, -2.5) circle (0.25 cm);
\node at (0.5, -2.5) {$\mu_1$};
\fill[white] (-0.5, -3.5) circle (0.25 cm);
\node at (-0.5, -3.5) {$\mu_1$};

\fill[white] (2.5, 1.5) circle (0.25 cm);
\node at (2.5, 1.5) {$\mu_2$};
\fill[white] (1.5, 0.5) circle (0.25 cm);
\node at (1.5, 0.5) {$\mu_2$};
\fill[white] (0.5, -0.5) circle (0.25 cm);
\node at (0.5, -0.5) {$\mu_2$};
\fill[white] (-0.5, -1.5) circle (0.25 cm);
\node at (-0.5, -1.5) {$\mu_2$};
\fill[white] (-1.5, -2.5) circle (0.25 cm);
\node at (-1.5, -2.5) {$\mu_2$};

\fill[white] (1.5, 2.5) circle (0.25 cm);
\node at (1.5, 2.5) {$\mu_3$};
\fill[white] (0.5, 1.5) circle (0.25 cm);
\node at (0.5, 1.5) {$\mu_3$};
\fill[white] (-0.5, 0.5) circle (0.25 cm);
\node at (-0.5, 0.5) {$\mu_3$};
\fill[white] (-1.5, -0.5) circle (0.25 cm);
\node at (-1.5, -0.5) {$\mu_3$};
\fill[white] (-2.5, -1.5) circle (0.25 cm);
\node at (-2.5, -1.5) {$\mu_3$};

\fill[white] (0.5, 3.5) circle (0.25 cm);
\node at (0.5, 3.5) {$\mu_4$};
\fill[white] (-0.5, 2.5) circle (0.25 cm);
\node at (0-.5, 2.5) {$\mu_4$};
\fill[white] (-1.5, 1.5) circle (0.25 cm);
\node at (-1.5, 1.5) {$\mu_4$};
\fill[white] (-2.5, 0.5) circle (0.25 cm);
\node at (-2.5, 0.5) {$\mu_4$};
\fill[white] (-3.5, -0.5) circle (0.25 cm);
\node at (-3.5, -0.5) {$\mu_4$};

\end{tikzpicture}} 

\caption{The restriction of the octahedron recurrence to the walls of $T_4$. 
The four corners of the figure are labeled by the vertices of the tetrahedron $T_4$. Fold along the dashed line; map the left side linearly to the wall $i=0$ and map the right side linearly to the wall $j=0$.
Each edge is labelled with $\tilh_{\text{left endpt}} - \tilh_{\text{right endpt}}$.
The red edges indicate  $\SHiveBottom$, $\SSkepBottom$,  $\SSkepTop$ and $\SHiveTop$.
}
\label{WallShapes}
\end{figure}

\begin{proof}[Proof of the easy parts of Theorem~\ref{TechnicalOctTheorem}]
The easy parts are that condition~(5) implies conditions~(1), (2), (3) and (4), and the claims about the boundary values.

We  check that $(5)\! \implies\! (1)$. The hive inequalities~\eqref{HiveInequalities} turn into the following conditions:
\[ \begin{array}{lcl}
\tilh_{(i+1)j(n-i-j-1)} + \tilh_{i(j+1)(n-i-j-1)} &\geq& \tilh_{ij(n-i-j)} + \tilh_{(i+1)(j+1)(n-i-j-2)} \\
\tilh_{i(j+1)(n-i-j-1)} + \tilh_{ij(n-i-j)} &\geq& \tilh_{(i-1)(j+1)(n-i-j)} + \tilh_{(i+1)j(n-i-j-1)} \\
\tilh_{(i+1)j(n-i-j-1)} + \tilh_{ij(n-i-j)} &\geq& \tilh_{(i+1)(j-1)(n-i-j)} + \tilh_{i(j+1)(n-i-j-1)} \\
\end{array}\]
All of these are rhombus inequalities.
The proof that $(5)\! \implies\! (2)$ is the same with $t$ negated. 

We now show that $(5)\! \implies \!(3)$.  
In each of the following inequalities, let $r$ be $0$ if $(i,j) \in \Delta^+_n$ and $1$ if $(i,j) \in \Delta^-_n$.
The first four skep inequalities~\eqref{SkepInequalities} turn into the following conditions:
\[ \begin{array}{lcl}
\tilh_{ijr} - \tilh_{(i+1)(j+1)r} &\leq& \tilh_{(i+1)j(1-r)} - \tilh_{(i+2)(j+1)(1-r)} \\
\tilh_{ijr} - \tilh_{(i+1)(j+1)r} &\leq& \tilh_{i(j+1)(1-r)} - \tilh_{(i+1)(j+2)(1-r)} \\
\tilh_{ijr} - \tilh_{(i+1)(j-1)r} &\leq& \tilh_{i(j-1)(1-r)} - \tilh_{(i+1)(j-2)(1-r)} \\
\tilh_{ij}r - \tilh_{(i+1)(j-1)r} &\leq& \tilh_{(i+1)j(1-r)} - \tilh_{(i+2)(j-1)(1-r)} \\
\end{array}\]
All of these are rhombus inequalities.

It remains to handle the inequality $0 \leq g_{00} - g_{11}$. Put $s = 0$ if $n$ is even and $1$ if $n$ is odd; then this inequality turns into the condition $0 \leq \tilh_{00s} - \tilh_{11s}$. We now use the final case of the octahedron recurrence~\eqref{OctRecurr} to write  $\tilh_{00s} = \tilh_{10(1-s)} + \tilh_{01(1-s)}- \tilh_{00(2-3s)}$, so we can rewrite our condition as
\[ 0 \leq \left( \tilh_{10(1-s)} + \tilh_{01(1-s)} -  \tilh_{00(2-3s)} \right)- \tilh_{11s} .\]
This is also a rhombus inequality.  The proof that $(5)\! \implies\! (4)$ is the same with $t$ negated. 

We next check the conditions on the boundary values.
First of all, the four arrays all have the same diagonal boundary, because all of the diagonal boundaries are
\[ ( \tilh_{(n-1)10} - \tilh_{n00},\  \tilh_{(n-2)20} - \tilh_{(n-1)10},\   \tilh_{(n-3)30} - \tilh_{(n-2)20},\  \ldots,\   \tilh_{0n0} - \tilh_{1(n-1)0}). \]

We now check the conditions on the side boundaries. Let
\[ B = \{ (i,0,t) : |t| \leq n-i,\ i+t \equiv n \bmod 2 \} \cup  \{ (0,j,t) : |t| \leq n-j,\ j+t \equiv n \bmod 2 \} .\]
$B$ is located on the two side walls of the tetrahedron $T_n$, and the various side boundaries are $\tilde{h}$ restricted to various subsets of $B$.
Figure~\ref{WallShapes} depicts $B$ for $n=4$, with the corners labeled by their coordinates in $T_n$. 
A dashed line shows $i=j=0$; we fold along this line to wrap our diagram onto the two walls of $T_n$.
Many of the differences between values of $\tilh|_B$ are equal to each other, as indicated by the edge labels; we verify this now.

Let 
\[ D = \{ (k,t) \in \ZZ^2 : |k| + |t| \leq n,\ k + t \equiv n \bmod 2 \} .\]
Define a bijection $\phi : D \to B$ by
\[ \phi(k,t) = \begin{cases} (k,0,t) & k \geq 0 \\ (0,-k,t) & k \leq 0 \end{cases} .\]
In otherwise, $\phi$ is the map which wraps the planar diagram in Figure~\ref{WallShapes} onto the two walls of $T_4$.
Let $d$ be the composition $D \overset{\phi}{\longrightarrow} B \overset{\tilh}{\longrightarrow} \ZZ$. 
Then the bottom three cases of the octahedron recurrence~\eqref{OctRecurr} all collapse to one equation:
\begin{equation} d_{kt} + d_{k(t-2)} = d_{(k-1)(t-1)} + d_{(k+1)(t-1)} . \label{dEquality} \end{equation}

Let $\HB(\tilh|_{\SHiveBottom}) = (\lambda_1, \lambda_2, \ldots, \lambda_n)$ and $\VB(\tilh|_{\SHiveBottom}) = (\mu_1, \mu_2, \ldots, \mu_n)$. 
Translated into $d$-values, this says that
\[ 
\begin{array}{lcl}
(d_{(n-1)(-1)} - d_{n0},\ d_{(n-2)(-2)} - d_{(n-1)(-1)},\ \ldots,\ d_{0(-n)} - d_{1(-n+1)}) &=& (\lambda_1,\ \ldots,\ \lambda_n) \\
(d_{(-1)(-n+1)}-d_{0(-n)},\ d_{(-2)(-n+2)}-d_{(-1)(-n)},\ \ldots,\ d_{(-n)0} - d_{(-n+1)(-1)})  &=& (\mu_1,\ \ldots,\ \mu_n) \\
\end{array}
\]
These are the edge labels on the bottom of Figure~\ref{WallShapes}.

We can rewrite~\eqref{dEquality} as $d_{(k-1)(t-1)} - d_{kt} =  d_{k(t-2)}-d_{(k+1)(t-1)}$. Thus, $d_{(k-1)(t-1)} - d_{kt}$ depends only on $k+t$, and must be $\lambda_{(n-k-t)/2+1}$
This formula is much clearer visually: In Figure~\ref{WallShapes}, edges in direction $(1,1)$ which lie on parallel sides of a square must be labeled with the same value of $\lambda$.
Similarly, $d_{(k-1)(t-1)} - d_{k(t-2)} =  d_{kt}-d_{(k+1)(t-1)}$ depends only on $k-t$, and must be $\mu_{(n-k+t)/2}$. In other words, edges in direction $(1,-1)$ which lie on parallel sides of a square must be labeled with the same value of $\mu$.

We now look at the various intersections $\SHiveBottom \cap B$, $\SSkepBottom \cap B$,  $\SSkepTop \cap B$ and $\SHiveTop \cap B$.
These are the red paths from $(n,0,0)$ to $(0,n,0)$ in Figure~\ref{WallShapes}, in the order from bottom to top of the figure.
 We see that the boundary values are permuted exactly as claimed.
\end{proof}

We now turn to the hard part of Theorem~\ref{TechnicalOctTheorem}; that any of (1), (2), (3) or (4) implies (5). 
This involves introducing more terminology from~\cite{HK}: \newword{sections} and \newword{wavefronts}.

\begin{definition}
Write \[ \RR \Delta_n:= \Hull((n,0), (0,n), (0,0))\  \text{and} \  \RR T_n := \Hull((n,0,0), (0,n,0), (0,0,n), (0,0,(-n))).\]
Recall the definition of a unit triangle from Definition~\ref{DefnTrianglesRhombiInequalities}. A \newword{section} $\cS$ is a union of unit triangles, inside $\RR T_n$, such that $\proj : \cS \to \RR \Delta_n$ is bijective. 
\end{definition}

Thus, if $\cS$ is a section, then $\cS \cap T_n \xrightarrow{\proj} \Delta_n$ is bijective. 
However, we emphasize that $\cS$ is \emph{not} determined by the discrete set $\cS \cap T_n$.
There are four sections which will play a key role in our proof of Theorem~\ref{TechnicalOctTheorem}:
\begin{definition}
First, define
\[ \cSTop = \Hull((n,0,0), (0,n,0), (0,0,n)) . \]
So $\cSTop \cap T_n = \SHiveTop$. We now define two sections with $\cS_1 \cap T_n = \cS_2 \cap T_n = \SSkepTop$.

Recall the function $\epsilon$ from Definition~\ref{KeySlices}.
We make the following definitions:
\[ \begin{array}{lcl}
\NWTri{i}{j} &=& \Hull((i,j,\epsilon(i,j)),\ (i+1,j,1-\epsilon(i,j)),\ (i,j-1,1-\epsilon(i,j))) \\
\NETri{i}{j} &=& \Hull((i,j,\epsilon(i,j)),\ (i-1,j,1-\epsilon(i,j)),\ (i,j-1,1-\epsilon(i,j))) \\
\SETri{i}{j} &=& \Hull((i,j,\epsilon(i,j)),\ (i-1,j,1-\epsilon(i,j)),\ (i,j+1,1-\epsilon(i,j))) \\
\SWTri{i}{j} &=& \Hull((i,j,\epsilon(i,j)),\ (i+1,j,1-\epsilon(i,j)),\ (i,j+1,1-\epsilon(i,j)))  \\
\end{array} \]
So the triangle depicts the projection to the $(i,j)$ plane, and the index gives the coordinates of the right angled vertex of the projection.
We now define:\pagebreak[2]
\[
\begin{array}{lcl}
\cS_1 &=& \displaystyle{ \bigcup_{i=0}^{n-1} \bigcup_{j=0}^{n-1-i}\  \SWTri{i}{j} \ \cup \ \bigcup_{i=1}^{n-1} \bigcup_{j=1}^{n-i} \NETri{i}{j} } \\[0.75 cm]
\cS_2 &=& \displaystyle{ \bigcup_{i=0}^{n-1} \bigcup_{j=1}^{n-1-i}\  \NWTri{i}{j} \ \cup \ \bigcup_{i=1}^{n-1} \bigcup_{j=0}^{n-1-i} \SETri{i}{j} } \ \cup \ \bigcup_{k=0}^{n-1} \ \SWTri{k}{(n-1-k)} \\
\end{array} \]

Finally, let $s = 0$ for $n$ even and $s=1$ for $n$ odd. We put
\[ \cS'_1 = \cS_1 \setminus  \SWTri{0}{0} \ \cup\  \Hull((0,0,2-3s),(1,0,1-s),(0,1,1-s)) . \]
In other words, $\cS'_1$ deletes $\SWTri{0}{0} = \Hull((0,0,s),(1,0,1-s),(0,1,1-s))$ and replaces it with another triangle with the same hypotenuse and same projection onto the $(i,j)$ plane.
\end{definition}

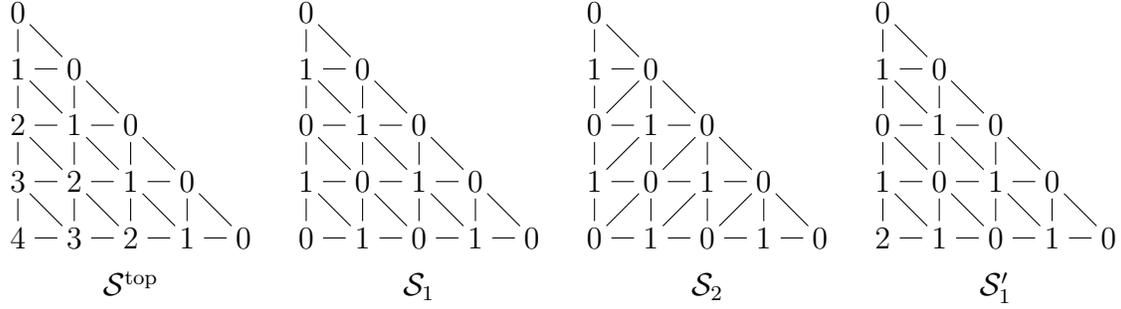
\begin{figure}
\centerline{
\begin{tabular}{cccc}
\begin{tikzpicture}[scale=0.75] 
\fill[white] (4,0) circle (0.3 cm);
\node at (4,0) {$0$};
\fill[white] (3,0) circle (0.3 cm); \fill[white] (3,1) circle (0.3 cm);
\node at (3,0) {$1$}; \node at (3,1) {$0$};
\fill[white] (2,0) circle (0.3 cm);\fill[white] (2,1) circle (0.3 cm); \fill[white] (2,2) circle (0.3 cm);
\node at (2,0) {$2$}; \node at (2,1) {$1$}; \node at (2,2) {$0$};
\fill[white] (1,0) circle (0.3 cm);\fill[white] (1,1) circle (0.3 cm); \fill[white] (1,2) circle (0.3 cm); \fill[white] (1,3) circle (0.3 cm);
\node at (1,0) {$3$}; \node at (1,1) {$2$}; \node at (1,2) {$1$}; \node at (1,3) {$0$};
\fill[white] (0,0) circle (0.3 cm);\fill[white] (0,1) circle (0.3 cm); \fill[white] (0,2) circle (0.3 cm); \fill[white] (0,3) circle (0.3 cm); \fill[white] (0,4) circle (0.3 cm);
\node at (0,0) {$4$}; \node at (0,1) {$3$}; \node at (0,2) {$2$}; \node at (0,3) {$1$}; \node at (0,4) {$0$};
\begin{pgfonlayer}{bg}    
\draw (0,0) -- (0,4);  \draw (1,0) -- (1,3); \draw (2,0) -- (2,2);   \draw (3,0) -- (3,1); 
\draw (0,0) -- (4,0);  \draw (0,1) -- (3,1); \draw (0,2) -- (2,2);   \draw (0,3) -- (1,3); 
\draw (4,0) -- (0,4); \draw (3,0) -- (0,3); \draw (2,0) -- (0,2); \draw (1,0) -- (0,1); 
\end{pgfonlayer}
\end{tikzpicture}  & 
\begin{tikzpicture}[scale=0.75] 
\fill[white] (4,0) circle (0.3 cm);
\node at (4,0) {$0$};
\fill[white] (3,0) circle (0.3 cm); \fill[white] (3,1) circle (0.3 cm);
\node at (3,0) {$1$}; \node at (3,1) {$0$};
\fill[white] (2,0) circle (0.3 cm);\fill[white] (2,1) circle (0.3 cm); \fill[white] (2,2) circle (0.3 cm);
\node at (2,0) {$0$}; \node at (2,1) {$1$}; \node at (2,2) {$0$};
\fill[white] (1,0) circle (0.3 cm);\fill[white] (1,1) circle (0.3 cm); \fill[white] (1,2) circle (0.3 cm); \fill[white] (1,3) circle (0.3 cm);
\node at (1,0) {$1$}; \node at (1,1) {$0$}; \node at (1,2) {$1$}; \node at (1,3) {$0$};
\fill[white] (0,0) circle (0.3 cm);\fill[white] (0,1) circle (0.3 cm); \fill[white] (0,2) circle (0.3 cm); \fill[white] (0,3) circle (0.3 cm); \fill[white] (0,4) circle (0.3 cm);
\node at (0,0) {$0$}; \node at (0,1) {$1$}; \node at (0,2) {$0$}; \node at (0,3) {$1$}; \node at (0,4) {$0$};
\begin{pgfonlayer}{bg}    
\draw (0,0) -- (0,4);  \draw (1,0) -- (1,3); \draw (2,0) -- (2,2);   \draw (3,0) -- (3,1); 
\draw (0,0) -- (4,0);  \draw (0,1) -- (3,1); \draw (0,2) -- (2,2);   \draw (0,3) -- (1,3); 
\draw (4,0) -- (0,4); \draw (3,0) -- (0,3); \draw (2,0) -- (0,2); \draw (1,0) -- (0,1); 
\end{pgfonlayer}
\end{tikzpicture}    & 
\begin{tikzpicture}[scale=0.75] 
\fill[white] (4,0) circle (0.3 cm);
\node at (4,0) {$0$};
\fill[white] (3,0) circle (0.3 cm); \fill[white] (3,1) circle (0.3 cm);
\node at (3,0) {$1$}; \node at (3,1) {$0$};
\fill[white] (2,0) circle (0.3 cm);\fill[white] (2,1) circle (0.3 cm); \fill[white] (2,2) circle (0.3 cm);
\node at (2,0) {$0$}; \node at (2,1) {$1$}; \node at (2,2) {$0$};
\fill[white] (1,0) circle (0.3 cm);\fill[white] (1,1) circle (0.3 cm); \fill[white] (1,2) circle (0.3 cm); \fill[white] (1,3) circle (0.3 cm);
\node at (1,0) {$1$}; \node at (1,1) {$0$}; \node at (1,2) {$1$}; \node at (1,3) {$0$};
\fill[white] (0,0) circle (0.3 cm);\fill[white] (0,1) circle (0.3 cm); \fill[white] (0,2) circle (0.3 cm); \fill[white] (0,3) circle (0.3 cm); \fill[white] (0,4) circle (0.3 cm);
\node at (0,0) {$0$}; \node at (0,1) {$1$}; \node at (0,2) {$0$}; \node at (0,3) {$1$}; \node at (0,4) {$0$};
\begin{pgfonlayer}{bg}    
\draw (0,0) -- (0,4);  \draw (1,0) -- (1,3); \draw (2,0) -- (2,2);   \draw (3,0) -- (3,1); 
\draw (0,0) -- (4,0);  \draw (0,1) -- (3,1); \draw (0,2) -- (2,2);   \draw (0,3) -- (1,3); 
\draw (4,0) -- (0,4); 
\draw (1,0) -- (2,1); \draw (2,0) -- (3,1);
\draw (0,0) -- (2,2);
\draw (0,1) -- (1,2); \draw (0,2) -- (1,3);
\end{pgfonlayer}
\end{tikzpicture}   & 
\begin{tikzpicture}[scale=0.75] 
\fill[white] (4,0) circle (0.3 cm);
\node at (4,0) {$0$};
\fill[white] (3,0) circle (0.3 cm); \fill[white] (3,1) circle (0.3 cm);
\node at (3,0) {$1$}; \node at (3,1) {$0$};
\fill[white] (2,0) circle (0.3 cm);\fill[white] (2,1) circle (0.3 cm); \fill[white] (2,2) circle (0.3 cm);
\node at (2,0) {$0$}; \node at (2,1) {$1$}; \node at (2,2) {$0$};
\fill[white] (1,0) circle (0.3 cm);\fill[white] (1,1) circle (0.3 cm); \fill[white] (1,2) circle (0.3 cm); \fill[white] (1,3) circle (0.3 cm);
\node at (1,0) {$1$}; \node at (1,1) {$0$}; \node at (1,2) {$1$}; \node at (1,3) {$0$};
\fill[white] (0,0) circle (0.3 cm);\fill[white] (0,1) circle (0.3 cm); \fill[white] (0,2) circle (0.3 cm); \fill[white] (0,3) circle (0.3 cm); \fill[white] (0,4) circle (0.3 cm);
\node at (0,0) {$2$}; \node at (0,1) {$1$}; \node at (0,2) {$0$}; \node at (0,3) {$1$}; \node at (0,4) {$0$};
\begin{pgfonlayer}{bg}    
\draw (0,0) -- (0,4);  \draw (1,0) -- (1,3); \draw (2,0) -- (2,2);   \draw (3,0) -- (3,1); 
\draw (0,0) -- (4,0);  \draw (0,1) -- (3,1); \draw (0,2) -- (2,2);   \draw (0,3) -- (1,3); 
\draw (4,0) -- (0,4); \draw (3,0) -- (0,3); \draw (2,0) -- (0,2); \draw (1,0) -- (0,1); 
\end{pgfonlayer}
\end{tikzpicture} \\
$\cSTop$ & $\cS_1$ & $\cS_2$ & $\cS'_1$ \\
\end{tabular}}
\caption{The projections of $\cSTop$, $\cS_1$, $\cS_2$ and $\cS'_1$ to the $(i,j)$-plane (for $n=4$). The numbers indicate the value of the $t$-coordinate.} \label{KeySections}
\end{figure}

Figure~\ref{KeySections} shows the projections of $\cSTop$, $\cS_1$, $\cS_2$ and $\cS'_1$ onto the $(i,j)$ plane, with the numbers indicating the values of the $t$-coordinate.

\begin{definition} \label{DefinitionWavefront}
A \newword{wavefront} is a two-dimensional subset of $\RR T_n$ of one of these  kinds:
\begin{enumerate}
\item The intersection of $\RR T_n$ with $\{ t-c = i-j \}$, for $c \equiv n \bmod 2$, $|c|<n$.
\item The intersection of $\RR T_n$ with  $\{t-c = -i+j\}$, for $c \equiv n \bmod 2$,  $|c|<n$.
\item The intersection of $\RR T_n$ with $\{ |t-c| = i+j \}$ for $c \equiv n \bmod 2$,  $|c|<n$. 
\end{enumerate}
In the last case, the point $(0,0,c)$ is called the \newword{cutpoint} of the wavefront. 
\end{definition}

\begin{remark} 
Henriques and Kamnitzer work with an octahedron recurrence which is defined in an infinite square prism, not just in the finite tetrahedron $\RR T_n$. In that context, the three cases in the definition of wavefronts collapse to one case. We have intersected the wavefronts with $\RR T_n$, at the expense of breaking symmetry.
\end{remark}

\begin{definition}
We say that a wavefront $\cW$ and a section $\cS$ are \newword{transverse} if
\begin{enumerate}
\item $\cW \cap S$ has dimension $\leq 1$ and
\item If $\cW$ is as in the last case of Definition~\ref{DefinitionWavefront}, then $\cS$ does not contain the cut point.
\end{enumerate}
Given $\tilh : T_n \to \ZZ$, a wavefront $\cW$ and a section $\cS$ transverse to $\cW$, we'll say that $\tilh$ \newword{satisfies the rhombus inequalities along $\cW \cap \cS$} if $\tilh$ satisfies the rhombus inequalities for those rhombi $R$ such that 
\begin{enumerate} \item $R \subset \cS$ and \item the short diagonal of $R$ is contained in $\cW \cap \cS$. \end{enumerate}.
\end{definition}

\begin{lemma} \label{EssentiallyHK}
Let $\tilh : T_n \to \ZZ$ be a function obeying the octahedron recurrence. Suppose that, for every wavefront $\cW$, there is a section $\cS(\cW)$ transverse to $\cW$ such that  $\tilh$ satisfies rhombus inequalities along $\cW \cap \cS(\cW)$. Then $\tilh$ satisfies all rhombus inequalities.
\end{lemma}

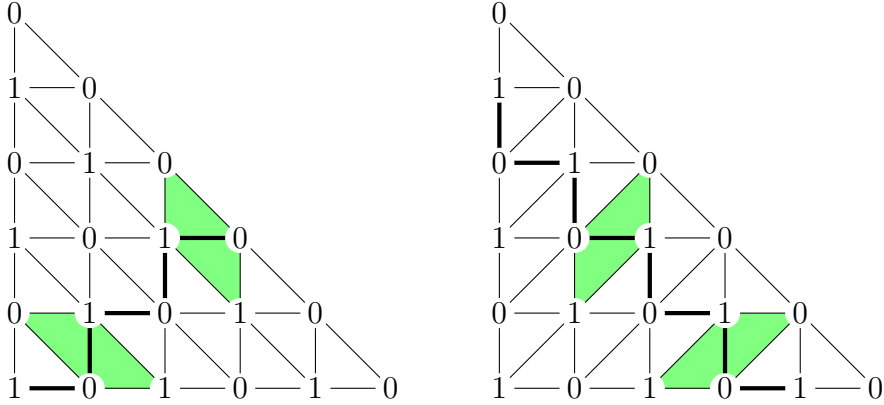
\begin{figure}
\centerline{
\begin{tikzpicture} 
\fill[white] (5,0) circle (0.2 cm);
\node at (5,0) {$0$};
\fill[white] (4,0) circle (0.2 cm); \fill[white] (4,1) circle (0.2 cm);
\node at (4,0) {$1$}; \node at (4,1) {$0$};
\fill[white] (3,0) circle (0.2 cm);\fill[white] (3,1) circle (0.2 cm); \fill[white] (3,2) circle (0.2 cm);
\node at (3,0) {$0$}; \node at (3,1) {$1$}; \node at (3,2) {$0$};
\fill[white] (2,0) circle (0.2 cm);\fill[white] (2,1) circle (0.2 cm); \fill[white] (2,2) circle (0.2 cm); \fill[white] (2,3) circle (0.2 cm);
\node at (2,0) {$1$}; \node at (2,1) {$0$}; \node at (2,2) {$1$}; \node at (2,3) {$0$};
\fill[white] (1,0) circle (0.2 cm);\fill[white] (1,1) circle (0.2 cm); \fill[white] (1,2) circle (0.2 cm); \fill[white] (1,3) circle (0.2 cm); \fill[white] (1,4) circle (0.2 cm);
\node at (1,0) {$0$}; \node at (1,1) {$1$}; \node at (1,2) {$0$}; \node at (1,3) {$1$}; \node at (1,4) {$0$};
\fill[white] (0,0) circle (0.2 cm);\fill[white] (0,1) circle (0.2 cm); \fill[white] (0,2) circle (0.2 cm); \fill[white] (0,3) circle (0.2 cm); \fill[white] (0,4) circle (0.2 cm); \fill[white] (0,5) circle (0.2 cm);
\node at (0,0) {$1$}; \node at (0,1) {$0$}; \node at (0,2) {$1$}; \node at (0,3) {$0$}; \node at (0,4) {$1$}; \node at (0,5) {$0$};
\begin{pgfonlayer}{bg}    
     \path [fill=green!50] (0,1) -- (1,0) -- (2,0) -- (1,1);
          \path [fill=green!50] (3,1) -- (2,2) -- (2,3) -- (3,2);
\draw (0,0) -- (0,5);  \draw (1,0) -- (1,4); \draw (2,0) -- (2,3);   \draw (3,0) -- (3,2);   \draw (4,0) -- (4,1); 
\draw (0,0) -- (5,0);  \draw (0,1) -- (4,1); \draw (0,2) -- (3,2);   \draw (0,3) -- (2,3);   \draw (0,4) -- (1,4); 
\draw (5,0) -- (0,5); \draw (4,0) -- (0,4);  \draw (3,0) -- (0,3); \draw (2,0) -- (0,2); \draw (1,0) -- (0,1); 
\draw[ultra thick] (0,0) -- (1,0) -- (1,1) -- (2,1) -- (2,2) -- (3,2);
\end{pgfonlayer}
\end{tikzpicture}   \qquad 
\begin{tikzpicture} 
\fill[white] (5,0) circle (0.2 cm);
\node at (5,0) {$0$};
\fill[white] (4,0) circle (0.2 cm); \fill[white] (4,1) circle (0.2 cm);
\node at (4,0) {$1$}; \node at (4,1) {$0$};
\fill[white] (3,0) circle (0.2 cm);\fill[white] (3,1) circle (0.2 cm); \fill[white] (3,2) circle (0.2 cm);
\node at (3,0) {$0$}; \node at (3,1) {$1$}; \node at (3,2) {$0$};
\fill[white] (2,0) circle (0.2 cm);\fill[white] (2,1) circle (0.2 cm); \fill[white] (2,2) circle (0.2 cm); \fill[white] (2,3) circle (0.2 cm);
\node at (2,0) {$1$}; \node at (2,1) {$0$}; \node at (2,2) {$1$}; \node at (2,3) {$0$};
\fill[white] (1,0) circle (0.2 cm);\fill[white] (1,1) circle (0.2 cm); \fill[white] (1,2) circle (0.2 cm); \fill[white] (1,3) circle (0.2 cm); \fill[white] (1,4) circle (0.2 cm);
\node at (1,0) {$0$}; \node at (1,1) {$1$}; \node at (1,2) {$0$}; \node at (1,3) {$1$}; \node at (1,4) {$0$};
\fill[white] (0,0) circle (0.2 cm);\fill[white] (0,1) circle (0.2 cm); \fill[white] (0,2) circle (0.2 cm); \fill[white] (0,3) circle (0.2 cm); \fill[white] (0,4) circle (0.2 cm); \fill[white] (0,5) circle (0.2 cm);
\node at (0,0) {$1$}; \node at (0,1) {$0$}; \node at (0,2) {$1$}; \node at (0,3) {$0$}; \node at (0,4) {$1$}; \node at (0,5) {$0$};
\begin{pgfonlayer}{bg}    
          \path [fill=green!50] (2,0) -- (3,0) -- (4,1) -- (3,1) -- (2,0);
          \path [fill=green!50] (1,1) -- (2,2) -- (2,3) -- (1,2) -- (1,1);
\draw (0,0) -- (0,5);  \draw (1,0) -- (1,4); \draw (2,0) -- (2,3);   \draw (3,0) -- (3,2);   \draw (4,0) -- (4,1); 
\draw (0,0) -- (5,0);  \draw (0,1) -- (4,1); \draw (0,2) -- (3,2);   \draw (0,3) -- (2,3);   \draw (0,4) -- (1,4); 
\draw (5,0) -- (0,5); 
\draw (3,0) -- (4,1);  \draw (2,0) -- (3,1); \draw (1,0) -- (3,2); \draw (0,0) -- (2,2);  \draw (0,1) -- (2,3); \draw (0,2) -- (1,3);  \draw (0,3) -- (1,4);
\draw[ultra thick] (4,0) -- (3,0) -- (3,1) -- (2,1) -- (2,2) -- (1,2)-- (1,3) -- (0,3) -- (0,4);
\end{pgfonlayer}
\end{tikzpicture}  
}
\caption{The intersections $\cW_1 \cap \cS_1$ and $\cW_2 \cap \cS_2$ (in bold), and some characteristic rhombus inequalities (in green). In the figure, we have $n=5$, $\cW_1 = \{ t-1 = -i+j \}$ and $\cW_2 = \{ |t+3| = i+j \}$.} \label{Intersections}
\end{figure}

\begin{proof}
This is morally~\cite[Lemma 3.1]{HK}. We show how to deduce this lemma from that one.

Henriques and Kamnitzer say that, if $\cS$ and $\cS'$ are two sections, with $\cS'$ ``in the future'' of $\cS$, and both $\cS$ and $\cS'$ are transverse to $\cW$, and if  $\tilh$ satisfies rhombus inequalities along $\cW \cap \cS$, then  $\tilh$ satisfies rhombus inequalities along $\cW \cap \cS'$.  However, we can discard the condition that $\cS'$ is in the future of $\cS$ because, for any sections $\cS$ and $\cS'$, the section $\cSTop$ is in the future of both $\cS$ and $\cS'$, and is transverse to every wavefront. So we can simply say that, if  $\cS$ and $\cS'$ are two sections which are transverse to $\cW$ and $\tilh$ satisfies rhombus inequalities along $\cW \cap \cS$, then  $\tilh$ satisfies rhombus inequalities along $\cW \cap \cS'$. 

Now, for any unit rhombus $R$, we can find a section $\cS_0$ containing $R$ and a wavefront $\cW$ transverse to $R$, with the short diagonal of $R$ contained in $\cW \cap \cS$.
So, if we can find a section $\cS(\cW)$ as above then we deduce that $\tilh$ satisfies rhombus inequalities along $\cW \cap \cS_0$, and thus in particular satisfies the rhombus inequality at $R$. 
So, if we can find such sections $\cS(\cW)$ for all wavefronts $\cW$, then $\tilh$ satisfies all rhombus inequalities.
\end{proof}

\begin{proof}[Proof of the remainder of Theorem~\ref{TechnicalOctTheorem}]
It remains to show that conditions (1) through (4) imply condition (5), the truth of all rhombus inequalities.

We first do the proof that $(1)\! \implies\! (5)$, although this is essentially in~\cite{HK}, since it will serve as a good warm up.
So, assume that $\tilh|_{\SHiveTop}$ is a hive. For every wavefront $\cW$, the section $\cSTop$ is transverse to $\cW$. 
And every rhombus inequality in $\cSTop$ is a hive inequality on $\tilh|_{\cSTop}$, so $\tilh$ satisfies the rhombus inequalities along $\cW \cap \cSTop$.
So, by Lemma~\ref{EssentiallyHK}, $\tilh$ satisfies all rhombus inequalities.
The proof that  $(2)\! \implies\! (5)$ is the same, negating all $t$-coordinates.

We now prove   $(3)\! \implies\! (5)$. There are three different kinds of wavefronts we must consider.

\textbf{Case 1:} Let $\cW$ be a wavefront of the form $t = c+i-j$ or $t = c-i+j$. Then $\cW$ is transverse to $\cS_1$. 
All unit rhombi contained in $\cS_1$ correspond to skep inequalities so, in particular, $\tilh$ satisfies the rhombus inequalities along $\cW \cap \cS_1$.

The left hand side of Figure~\ref{Intersections} shows an example of $\cW \cap \cS_1$ with bold lines. Two of the rhombus inequalities along $\cW \cap \cS_1$ are highlighted in green. 

\textbf{Case 2:}  Let $\cW$ be a wavefront of the form $|t-c| = i+j$ for $c \neq \epsilon(0,0)$. Then $\cW$ is transverse to $\cS_2$. 
Again, all unit rhombi contained in $\cS_2$ correspond to skep inequalities so, in particular, $\tilh$ satisfies the rhombus inequalities along $\cW \cap \cS_2$.

The right hand side of Figure~\ref{Intersections} shows $\cW \cap \cS_2$ with bold lines, where $\cW = \{ |t+3| = i+j \}$. (Note that $\{ t+3 = -(i+j) \}$ is disjoint from $\cS_2$, although it does intersect $\RR T_5$.)  Again, two of the rhombus inequalities along $\cW \cap \cS_2$ are highlighted in green. 

\textbf{Case 3:} Let $\cW$ be the wavefront $|t-\epsilon(0,0)| = i+j$. Abbreviate $\epsilon(0,0)$ to $s$.

In this case, $(0,0,s)$ is the cutpoint, so our section must not contain $(0,0,s)$. In particular, we cannot use a section whose vertices are $\SSkepTop$ this time. Instead, we use $\cS'_1$.

The intersection $\cW \cap \cS_1$ is the single line segment from $(1,0,1-s)$ to $(0,1,1-s)$. This is the short diagonal of the rhombus whose long diagonal runs from $(0,0,2-3s)$ to $(1,1,s)$, so we must verify that $\tilh(1,0,1-s) + \tilh(0,1,1-s) \geq \tilh(0,0,2-3s) + \tilh(1,1,s)$.

The last case of~\ref{OctRecurr} gives $\tilh(0,0,2-3s) =  \tilh(1,0,1-s) + \tilh(0,1,1-s) - \tilh(0,0,s)$
so we need to show that $0 \geq \tilh(1,1,s)  - \tilh(0,0,s)$.
This is the last skep inequality. 

The implication $(4)\! \implies\! (5)$ is the same with all of the $t$-coordinates negated. This concludes the proof of Theorem~\ref{TechnicalOctTheorem}.
\end{proof}

We have now proved Theorem~\ref{TechnicalOctTheorem}, so we have proved its consequences, Theorems~\ref{SkepLR}, \ref{SkepCommutative} and~\ref{SkepLRAsSum}.
We are now done with hives and with other connections to classical Littlewood-Richardson combinatorics. 
The rest of the paper works with skeps and  $L$-convexity.

\section{More on $L$-convexity} \label{LConvex2}

The goal of this section is to give alternate formulations of $L$-convexity.
To help the reader, we have included short proofs of known results.

\begin{definition}
Let $x$ and $y \in \ZZ^n$. We define
\[ 
\begin{array}{rcl}
\floorAve{x}{y} &:=& \left( \lfloor \tfrac{x_1+y_1}{2} \rfloor, \lfloor \tfrac{x_2+y_2}{2} \rfloor, \dots, \lfloor \tfrac{x_n+y_n}{2} \rfloor \right) \ \text{and} \\
\ceilAve{x}{y} &:=& \left( \lceil \tfrac{x_1+y_1}{2} \rceil, \lceil \tfrac{x_2+y_2}{2} \rceil, \dots, \lceil \tfrac{x_n+y_n}{2} \rceil \right) \\
x \inmin y &:=&  (\min(x_1, y_1), \min(x_2, y_2), \ldots, \min(x_n, y_n)) \ \text{and} \\ 
 x \inmax y &:=& (\max(x_1, y_1), \max(x_2, y_2), \ldots, \max(x_n, y_n))  .
 \end{array}
\]
Here $\lfloor z \rfloor$ is $z$ rounded down to an integer, and $\lceil z \rceil$ is $z$ rounded up to an integer.
\end{definition}
We note that $\floorAve{x}{y}$ $\ceilAve{x}{y}$, $x \inmin y$ and $x \inmax y$ are all in $\Pi(x,y)$, and that 
\[  \floorAve{x}{y}+\ceilAve{x}{y} = x \inmin y + x \inmax y = x+y . \]

\begin{definition}
Let  $e_1$, $e_2$, \ldots, $e_n$ be the standard basis of $\ZZ^n$.
For any subset $S$ of $[n]$, we define:
\[ e_S := \sum_{i \in S} e_i . \]
\end{definition}

\begin{theorem} \label{EquivalentLConvexSetDefinitions}
Let $K$ be a subset of $\ZZ^n$ such that $K+\One_n=K$. The following are equivalent.
\begin{enumerate}
\setcounter{enumi}{-1} 
\item There are values $c_{ij} \in \ZZ \cup \{\infty \}$, for $1 \leq i,j \leq n$, such that
\[ K = \{ x \in \ZZ^n : x_i - x_j \leq c_{ij} \} . \]
\item  For any $x$, $y \in K$, we have $\Pi(x,y) \subseteq K$.
\item For any $x$, $y \in K$, the vectors $x \inmin y$ and $x \inmax y$ are in $K$.
\item For any $x$, $y \in K$, the vectors  $\floorAve{x}{y}$ and $\ceilAve{x}{y}$ are in $K$.
\end{enumerate}
\end{theorem}

\begin{theorem} \label{EquivalentLConvexFunctionDefinitions}
Let $f : \ZZ^n \to \RR \cup \{ \infty \}$ be a function such that there exists a constant $r \in \RR$ with $f(x + \One_n) = f(x)+r$. The following are equivalent:
\begin{enumerate}
\item   For any $x$, $y$, $x'$, $y' \in \ZZ^n$ with $x+y=x'+y'$ and $x' \in \Pi(x,y)$, we have $f(x)+f(y) \geq f(x') + f(y')$.
\item For any $x$, $y \in \ZZ^n$, we have $f(x)+f(y) \geq f(x \inmin y)+f(x \inmax y)$.
\item For any $x$, $y \in \ZZ^n$, we have $f(x)+f(y) \geq f(\floorAve{x}{y})+f(\ceilAve{x}{y})$.
\item The set $\Finite(f):= \{ x\in \ZZ^n : f(x) < \infty \}$ is $L$-convex and, for any $x \in \ZZ^n$ and any subsets $I \subseteq J$ of $[n]$, we have $f(x) + f(x+e_I+e_J) \geq f(x+e_I) + f(x+e_J)$.
\end{enumerate}
\end{theorem}

\begin{definition}
A subset $K$ of  $\ZZ^n$ is called \newword{$L$-convex} if it obeys the equivalent conditions of Theorem~\ref{EquivalentLConvexSetDefinitions}; 
a function $f : \ZZ^n \to \RR \cup \{ \infty \}$ is called \newword{$L$-convex} if it obeys the equivalent conditions of Theorem~\ref{EquivalentLConvexFunctionDefinitions}.
\end{definition}

The following is obvious from every one of these perspectives:
\begin{cor} \label{Translate}
Let $K$ be a subset of $\ZZ^n$ and let $v$ be a vector in $\ZZ^n$. Then $K$ is $L$-convex if and only if $K+v$ is $L$-convex.
Likewise, let $f : \ZZ^n \to  \RR \cup \{ \infty \}$ be a function. Then $f$ is $L$-convex if and only if $x \mapsto f(x+v)$ is $L$-convex.
\end{cor}

\begin{remark} \label{Progression}
Condition~(3) of Theorem~\ref{EquivalentLConvexFunctionDefinitions} implies, in particular, that $f(x+z) + f(x-z) \geq 2 f(x)$ and, thus, the restriction of $f$ to any arithmetic progression in $\ZZ^n$ is convex. 
\end{remark}

\begin{remark}
Most of these conditions already appear in works of Murota.
Condition~(0) is discussed in \cite[Section 5.3]{Murota}.
Conditions~(1) are what we used in Definitions~\ref{LConvexSetDefn} and~\ref{LConvexFunctionDefn}. These conditions are the most natural for us, since they match the form of the LPP conjecture.
Conditions~(2) are used to define $L$-convexity in~\cite{Murota}, where they are called the $SBS[\ZZ]$ and $SBF[\ZZ]$ conditions.
These conditions  will be key to our proof of Theorem~\ref{PointCountLConcave}. 
Conditions~(3) are discussed in~\cite[Section 7.2]{Murota} and in~\cite{Midpoint}, where they are called ``discrete midpoint convexity". We do not use them, but we include them for comparison with the literature, and for the discussion in Remark~\ref{SchurLlogConcaveRemark}.
Condition~(4) is a special case of the parallelogram conditions from~\cite{Midpoint}, and also a special case of the $L^{\natural}-APR[\ZZ]$ condition from~\cite[Section 7.2]{Murota}. 
The importance of these particular parallelograms is referenced in~\cite[Proposition 7.5]{Murota}, which essentially states that~(4) implies~(2). 
\end{remark}

\begin{remark}
Of course, all of the criteria in Theorem~\ref{EquivalentLConvexFunctionDefinitions} turn into criteria for a function to be $L$-concave or to be $\Llog$-concave, simply replacing $f$ with $-f$ or with $- \log f$.
\end{remark}

\begin{remark}
The convex hulls of $L$-convex sets are called ``alcoved polytopes" in the work of Lam and Postnikov~\cite{LamPost1, LamPost2}. 
More specifically, these are the type $A$ alcoved polytopes, and there is a notion of ``alcoved polytope" for every root system.
The author thinks it would be valuable to develop the theory of $L$-convexity for the other Dynkin types.

\begin{remark}
In probability, $\Llog$-concave functions are called ``$\log$ supermodular functions" or are said to obey the ``strong FKG condition". 
See~\cite{FKG} or \cite[Chapter 6]{AlonSpencer}.
\end{remark}

\end{remark}

Before proving Theorems~\ref{EquivalentLConvexSetDefinitions} and ~\ref{EquivalentLConvexFunctionDefinitions}, we give a more explicit description of the geometry of $\Pi(x,y)$, and use it to introduce some notation:

\begin{lemma} \label{PiGeometry}
Let $x$ and $y \in \ZZ^n$. Then there is a chain of subsets  $[n] \supsetneq I_1 \supsetneq I_2 \supsetneq \cdots \supsetneq I_d \supsetneq \emptyset$,  positive integers $\ell_1$, $\ell_2$, \dots, $\ell_d$, and some $c \in \ZZ$, such that 
\[ y = x + \sum \ell_j e_{I_j} + c \One_n. \]
We then have
\[ \Pi(x,y) = x + [0,\ell_1] e_{I_1} +  [0,\ell_2] e_{I_2} + \cdots +  [0,\ell_d] e_{I_d} + \ZZ \One_n . \]
\end{lemma}

\begin{definition}
We set $L(x,y) := \sum \ell_i = \max_i(x_i - y_i) - \min_j(x_j - y_j)$ and call this the \newword{L-distance} from $x$ to $y$. 
L-distance is a pseudo-metric, meaning that $L(x,x)=0$, $L(x,y) = L(y,x)$ and $L(x,z) \leq L(x,y) + L(y,z)$, but $L(x,y) = 0$ doesn't imply that $x=y$, rather, it is equivalent to  $x - y \in \ZZ \One_n$.
We call $(I_1, \ldots, I_d)$ the \newword{directions from $x$ to $y$} and $(\ell_1, \ldots, \ell_n)$ the \newword{distances from $x$ to $y$}.
\end{definition}

\begin{remark}
Switching the roles of $x$ and $y$ replaces $(I_1, \ldots, I_d)$ with $([n] \setminus I_d, \ldots, [n] \setminus I_1)$ and replaces $(\ell_1, \ell_2, \ldots, \ell_d)$ with $(\ell_d, \ell_{d-1}, \ldots, \ell_1)$.
\end{remark}

\begin{remark} \label{PiGeometryRemark}
Geometrically, Lemma~\ref{PiGeometry} says that $\Pi(x,y)/\ZZ \One_n$ is the lattice points of a parallelepiped inside $\RR^n/\RR \One_n$. Specifically, the sides of this parallelepiped point in directions $e_{I_1}$, $e_{I_2}$, \dots, $e_{I_d}$ and have lengths $\ell_1$, $\ell_2$, \ldots, $\ell_d$. Figures~\ref{NotStrong} and~\ref{Example0358Figure} depict $\Pi(0000,0235)/\ZZ \One_4$ and $\Pi(0000, 0358)/\ZZ \One_4$. In both cases, the directions are $(\{ 2,3,4 \}, \{ 3,4 \},\ \{ 4 \})$; the distances are $(3,2,3)$ and $(2,1,2)$ respectively.

We have abbreviated $(0,2,3,5)$ to $0235$ above, etcetera. We will often shorten vectors in this way to improve readability in examples.
\end{remark}


\begin{proof}[Proof Of Lemma~\ref{PiGeometry}]
We start with preliminary reductions. We can translate by $x$, and hence we can assume that $x$ is $0$.
We can also reorder our coordinates such that $y_1 \leq y_2 \leq \cdots \leq y_n$. 
Let $h_0 < h_1 < \cdots < h_d$ be the set of distinct values of $y_i$. Let $I_j = \{ i : y_i \geq h_j \}$, including $I_0 = [n]$. We have
\[ y = h_0 \One_n + (h_1 - h_0) e_{I_1} + (h_2-h_1) e_{I_2} + \cdots + (h_d - h_{d-1}) e_{I_d}. \]
So, taking $c=h_0$ and $\ell_j = h_j - h_{j-1}$, we have the claimed formula for $y$.

We now check the claimed formula for $\Pi(x,y)$. We continue with the normalizations and definitions above, so $x=0$ and $y_1 \leq y_2 \leq \cdots \leq y_n$. 
It will be convenient to define $S_j = I_{j-1} \setminus I_j$. 
First, we check that $\Pi(x,y)$ is contained in the vector space spanned by the $e_{I_j}$, including $e_{I_0}$. Equivalently, we want to show that $\Pi(x,y)$ is in  the vector space spanned by the $e_{S_j}$. 
Indeed, let $a$ and $b$ in $S_j$. Then $x_a-x_b=y_a-y_b = 0$ so, if $z \in \Pi(x,y)$, then $0 \leq z_a - z_b \leq 0$ so $z_a = z_b$. So $z$ is constant on each block $S_j$, and  $\Pi(x,y)$ is in  the vector space spanned by the $e_{S_j}$ and, equivalently, in the span of the $e_{I_j}$.

Now, let $z = \sum c_j e_{I_j}$. We want to work out a criterion for when $z \in \Pi(x,y)$. The inequalities defining $\Pi(x,y)$ are $0 \leq z_b - z_a \leq y_b - y_a$ for $a \leq b$. In particular, if $a \in S_{j-1}$ and $b \in S_j$, then this becomes $0 \leq c_j  \leq \ell_j$. Conversely, suppose that $0 \leq c_j \leq \ell_j$ for $1 \leq j \leq d$. Let $1 \leq a \leq b \leq n$, with $a \in S_i$ and $b \in S_j$. Then $z_b - z_a = \sum_{k=i+1}^j c_k$ so $0 \leq z_b - z_a \leq \sum_{k=i+1}^j \ell_k = y_b - y_a$, so we have checked that the inequalities $0 \leq c_j \leq \ell_j$ imply all of the inequalities defining $\Pi(x,y)$. 
\end{proof}


\begin{lemma} \label{LHalfDistance}
For any $x$, $y \in \ZZ^n$, the L-distance $L(x, \floorAve{x}{y} )$ is either  $\lfloor \tfrac{L(x,y)}{2} \rfloor$ or $\lceil \tfrac{L(x,y)}{2} \rceil$. 
\end{lemma} 

\begin{proof}
Let $d = y-x$ and let $s$ be the vector obtained by sorting the entries of $d$ into decreasing order.
Let $z = \floorAve{x}{y}$, let $e = z-x$ and let $t$  be the vector obtained by sorting the entries of $e$ into decreasing order.
Then $e = \lfloor \tfrac{d}{2} \rfloor$ and  $t = \lfloor \tfrac{s}{2} \rfloor$. 
We have $L(x,y) = s_1 - s_n$ and $L(x,z) = t_1 - t_n =  \lfloor \tfrac{s_1}{2} \rfloor - \lfloor \tfrac{s_n}{2} \rfloor$

Running through the $4$ possible cases for $(s_1, s_2)$ modulo $2$, we see that $\lfloor \tfrac{s_1}{2} \rfloor - \lfloor \tfrac{s_n}{2} \rfloor$ is always either $\lfloor \tfrac{s_1 - s_n}{2} \rfloor$ or $\lceil \tfrac{s_1 - s_n}{2} \rceil$.
\end{proof}

\begin{proof}[Proof of Theorem~\ref{EquivalentLConvexSetDefinitions}]
The implication $(0) \!\implies\! (1)$ is immediate.
The implications $(1) \!\implies\! (2)$ and $(1) \!\implies\! (3)$ are immediate consequences of the observations that $x \inmin y$, $x \inmax y$, $\lfloor \tfrac{x+y}{2} \rfloor$ and $\lceil \tfrac{x+y}{2} \rceil$ are in $\Pi(x,y)$.
We will show $(2) \!\implies\! (0)$ and $(3) \!\implies\! (1)$, finishing the proof.

$\boldsymbol{(2) \!\implies\! (0)}$: This argument is taken from \cite[Proposition 5.2.(2)]{Murota}.

If $K = \emptyset$, just choose $c_{ij}$ such that $c_{12} + c_{21} < 0$. Assume from now on that $K \neq \emptyset$. Define $c_{ij} := \text{sup}_{q \in K} (q_i - q_j)$. Clearly, $K \subseteq \{ x \in \ZZ^n : x_i - x_j \leq c_{ij} \}$. We must show that, given any $p \in \{ x \in \ZZ^n : x_i - x_j \leq c_{ij} \}$, we have $p \in K$.

For each $1 \leq i,j \leq n$, we have $p_i - p_j \leq c_{ij}$, so there is some $q^{ij} \in K$ with $p_i - p_j \leq q^{ij}_i - q^{ij}_j$. Translating by $\One_n$, we may assume that $q^{ij}_i = p_i$, so $q^{ij}_j \leq p_j$.
Put $r^i = q^{i1} \inmin q^{i2} \inmin \cdots \inmin q^{in}$, so $r^i \in K$. We have $r^i_i = p_i \inmin p_i \inmin \cdots \inmin p_i = p_i$ and $r^i_j \leq q^{ij}_j \leq p_j$. Thus we have $r_1 \inmax r_2 \inmax \cdots \inmax r_n = p$ and thus $p \in K$.

$\boldsymbol{(3) \!\implies\! (1)}$: Our proof is by induction on $L(x,y)$. If $L(x,y)=0$, then $x-y \in \One_n$, $\Pi(x,y) = y+\ZZ \One_n$ and the result is clear.
Let $L(x,y) \geq 1$.

Define sequences of vectors  $x^0$, $x^1$, $x^2$,  \dots and $y^0$, $y^1$, $y^2$, \dots, by $(x^0, y^0) = (y,x)$ and $(x^{i+1}, y^{i+1}) = (\floorAve{x^i}{x}, \ceilAve{y^i}{y})$. By condition~(4), all of the $x^i$ and $y^i$ are in $K$. Also, it is easy to verify that $x^i+y^i = x+y$ for all $i$. Iterating Lemma~\ref{LHalfDistance} we have $L(x, x^{i+1}) = \lfloor \tfrac{L(x,x^i)}{2} \rfloor$ or $L(x, x^{i+1}) = \lceil \tfrac{L(x,x^i)}{2} \rceil$. Therefore, there must be some index $M$ for which $L(x, x^M) = 1$.

Since $L(x, x^M)=1$, we must have $x^M = x + e_J + c \One_n$ for some proper nonempty subset $J$ of $[n]$ and some $c \in \ZZ$. Since $x^M+y^M = x+y$, we also have $y^M = y-e_J - c \One_n$. Since $x^M \in \Pi(x,y)$, the set $J$ must be one of the directions from $x$ to $y$. Then $\Pi(x,y) = \Pi(x^M, y) \cup \Pi(x, y^M)$ and $L(x^M,y) = L(x, y^M) = L(x,y) -1$. By induction, we have $\Pi(x^M, y) \subseteq K$ and $\Pi(x, y^M) \subseteq K$, so $\Pi(x,y) =  \Pi(x^M, y) \cup \Pi(x, y^M) \subseteq K$ as desired.
\end{proof}

We now prove Theorem~\ref{EquivalentLConvexFunctionDefinitions}. 
Many of these implications are already known but we repeat the proofs both because they are short and because we will need to discuss these proofs in Remark~\ref{SchurLlogConcaveRemark}, where we discuss the possibility for a theory of Schur $\Llog$-concavity.

\begin{proof}[Proof of Theorem~\ref{EquivalentLConvexFunctionDefinitions}]
$\boldsymbol{(1) \implies (2)}$: Take $x'=x \inmin y$ and $y'=x \inmax y$.

$\boldsymbol{(2) \implies (3)}$: Our proof is by induction on $L(x,y)$. If $L(x,y)=0$ then $y = x + c \One_n$ for some $c \in \ZZ$, and $\left( \floorAve{x}{y}, \ceilAve{x}{y} \right) = \left( x + \lfloor \tfrac{c}{2} \rfloor \One_n,  y- \lfloor \tfrac{c}{2} \rfloor \One_n \right)$. 
Similarly, if  $L(x,y) = 1$, then $y = x + e_J + c \One_n$ for some $c \in \ZZ$. If $c$ is even, then $\left( \floorAve{x}{y}, \ceilAve{x}{y} \right) = \left( x + \tfrac{c}{2}  \One_n, y-\tfrac{c}{2} \One_n \right)$ and, if $c$ is odd, then $\left( \floorAve{x}{y}, \ceilAve{x}{y} \right) = \left( y- \lceil \tfrac{c}{2} \rceil  \One_n,  x +  \lceil \tfrac{c}{2} \rceil \One_n \right)$. These are our base cases.

Now, assume that $L(x,y) \geq 2$, so there are indices $i$ and $j$ with $x_i - y_i \geq x_j - y_j + 2$. Set $c = x_i - y_i-1$. Replacing $x$ by $\bar{x} := x - c \One_n$, we have $\bar{x}_i>y_i$ and $\bar{x}_j < y_j$. 
Set $x' = \bar{x} \inmin y$ and $y' = \bar{x} \inmax y$. Then $L(x', y') < L(\bar{x}, y) = L(x,y)$ so, by induction, $f(x')+f(y') \geq f(\floorAve{x'}{y'}) +  f(\ceilAve{x'}{y'}) =  f(\floorAve{\bar{x}}{y}) +  f(\ceilAve{\bar{x}}{y})$.
Also, by hypothesis~(2), we have $f(\bar{x})+f(y) \geq f(x')+f(y')$. So $f(\bar{x})+f(y) \geq f(\floorAve{\bar{x}}{y}) +  f(\ceilAve{\bar{x}}{y})$. Subtracting appropriate multiples of $c \One_n$ from each argument, we draw the desired conclusion.

%

$\boldsymbol{(3) \implies (4)}$. If $x$ and $y$ are in $\Finite(f)$, then $(3)$ says that $f(\floorAve{x}{y}) + f(\ceilAve{x}{y}) \leq f(x) + f(y) < \infty$, so $\floorAve{x}{y}$ and $\ceilAve{x}{y}$  are also in $\Finite(f)$. Thus, $\Finite(f)$ is $L$-convex according to condition~(3) in Theorem~\ref{EquivalentLConvexSetDefinitions}. We have  $\floorAve{x}{(x+e_I+e_J)} = x+e_I$ and $\ceilAve{x}{(x+e_I+e_J)} = x+e_J$, so $(3)$ implies the inequality in $(4)$.

$\boldsymbol{(4) \implies (1)}$.  Let $y = x + \sum \ell_j e_{I_j} + c \One_n$ as in Lemma~\ref{PiGeometry}. Then $x' = x + \sum a_j e_{I_j} + d \One_n$ for $a_j \in [0, \ell_j]$ and $d \in \ZZ$.
Without loss of generality, let $c=d=0$; this shifts both sides of the inequality by $cr$.
So $x' =  x+\sum a_j e_{I_j}$, $y' = x+ \sum (\ell_j-a_j) e_{I_j}$ and $y =  x+ \sum \ell_j e_{I_j}$. 

Let $p_1 p_2 \cdots p_A$ be the word $1^{a_1} 2^{a_2} \cdots d^{a_d}$ in the alphabet $[d]$, similarly, let $q_1 q_2 \cdots q_B$ be the word $1^{\ell_1-a_1} 2^{\ell_2-a_2} \cdots d^{\ell_d-a_d}$.
Put $u_{rs} = x+ \sum_{i=1}^r e_{I_{p_i}} +  \sum_{j=1}^s e_{I_{q_j}}$. 
So $u_{00} = x$, $u_{A0} = x+ \sum_{i=1}^A e_{I_{p_i}}  = x+ \sum_{k=1}^d a_k e_{I_k} = x'$ and, similarly, $u_{0B} = y'$ and $u_{AB} = y$. 

Note that all of the $u_{rs}$ lie in $x + \sum_j [0, \ell_j] e_{I_j}$ and hence in $\Pi(x,y)$. Since $\Finite(f)$ is assumed $L$-convex, we have $f(u_{rs}) < \infty$ for all $(r,s)$. Also, we have $u_{(r+1)s} = u_{rs} + e_{I_{p_{r+1}}}$, $u_{r(s+1)} = u_{rs} + e_{I_{q_{s+1}}}$ and $u_{(r+1)(s+1)} = u_{rs} + e_{I_{p_{r+1}}} + e_{I_{q_{s+1}}}$. So condition~(4) shows that
\begin{equation}
f(u_{rs}) + f(u_{(r+1)(s+1)}) \geq f(u_{(r+1)s}) + f(u_{r(s+1)}). \label{KeyRhombusIneq}
\end{equation}
Summing~\eqref{KeyRhombusIneq} over  $0 \leq r \leq A-1$ and $0 \leq s \leq B-1$, we get 
\[
\sum_{r=0}^{A-1} \sum_{s=0}^{B-1} f(u_{rs}) +  \sum_{r=1}^{A} \sum_{s=1}^{B} f(u_{rs}) \geq
\sum_{r=1}^{A} \sum_{s=0}^{B-1} f(u_{rs}) +  \sum_{r=0}^{A-1} \sum_{s=1}^{B} f(u_{rs}) .
\]
Cancelling the terms which occur on both sides (which are all finite!), we get $f(u_{00}) + f(u_{AB}) \geq f(u_{A0}) + f(u_{0B})$ or, in other words, $f(x)+f(y) \geq f(x') + f(y')$, as desired.
\end{proof}

\begin{eg} \label{4Implies1}
We write out how the proof of $(4) \implies (1)$ works in a particular case.
 We show that, if $f$ obeys $(3)$, then $f(0,0,0,0) + f(0,2,3,5) \geq f(0,2,2,4)+f(0,0,1,1)$. 
Our proof considers the $10$ quantities in the diagram below:
\[ \begin{tikzpicture}[scale=0.8]
  \node at (0,0) {$f(0,0,0,0)$};
    \node at (3,0) {$f(0,1,1,1)$};
        \node at (6,0) {$f(0,2,2,2)$};       
         \node at (9,0) {$f(0,2,2,3)$};
                  \node at (12,0) {$f(0,2,2,4)$};
  \node at (1.5,1.5) {$f(0,0,1,1)$};
    \node at (4.5,1.5) {$f(0,1,2,2)$};
        \node at (7.5,1.5) {$f(0,2,3,3)$};       
         \node at (10.5,1.5) {$f(0,2,3,4)$};
                  \node at (13.5,1.5) {$f(0,2,3,5)$};
\draw (1.25,0) -- (1.75,0);
\draw (4.25,0) -- (4.75,0);                  
\draw (7.25,0) -- (7.75,0);
\draw (10.25,0) -- (10.75,0);
\draw (2.75,1.5) -- (3.25,1.5);
\draw (5.75,1.5) -- (6.25,1.5);                  
\draw (8.75,1.5) -- (9.25,1.5);
\draw (11.75,1.5) -- (12.25,1.5);
\draw (0.25,0.25) -- (1.25, 1.25); \draw (3.25,0.25) -- (4.25, 1.25);  \draw (6.25,0.25) -- (7.25, 1.25);  \draw (9.25,0.25) -- (10.25, 1.25); \draw (12.25,0.25) -- (13.25, 1.25);
\end{tikzpicture}\]
In each of the four parallelograms, the sum of the values at the acute angles is greater than the sum of the values at the obtuse angles. Canceling the six values in the middle, we derive the conclusion.
The left hand side of Figure~\ref{NotStrong} depicts how these lattice points sit in space.
\end{eg}

\usetikzlibrary{math}
\usetikzlibrary{shapes}

\begin{figure}[b]
\centerline{
\begin{tikzpicture}[xscale=0.7, yscale=0.7]
\begin{pgfonlayer}{bg}    
\draw[thick] (2,1.5) -- (3.5,5);
\draw[thick] (4.75,1.5) -- (6.25,5);
\draw[thick] (1,0.75) -- (2.5, 4.25);
  \foreach \x in {0,1,...,2}
      \foreach \z in {0,1,...,1}
          \draw[thick] (2.75*\x+1.5*\z, 3.5*\z) -- (2.75*\x+2+1.5*\z, 1.5+3.5*\z);
  \foreach \y in {0,1,...,2}
      \foreach \z in {0,1,...,1}
      {
          \draw[white, fill=white] (\y+1.5*\z-0.1, 0.75*\y+3.5*\z-0.1) rectangle (5.5+\y+1.5*\z+0.1, 0.75*\y+3.5*\z+0.1);
          \draw[thick] (\y+1.5*\z, 0.75*\y+3.5*\z) -- (5.5+\y+1.5*\z, 0.75*\y+3.5*\z);
}           
\draw[thick] (0,0) -- (1.5, 3.5);
\draw[thick] (7.5,1.5) -- (9, 5);
\path[fill=white] (5.4, -0.1) -- (5.6, -0.1) -- (7.1,3.6) -- (6.9,3.6) -- (5.4, -0.1);
\draw[thick] (5.5, 0) -- (7, 3.5);
\path[fill=white] (2.65, -0.1) -- (2.75, -0.1) -- (4.35,3.6) -- (4.15, 3.6) -- (2.65,  -0.1);
\draw[thick] (2.75, 0) -- (4.25, 3.5);
\path[fill=white] (6.4, 0.74) -- (6.6, 0.74) -- (8.1,4.35) -- (7.9,4.35) -- (6.4, 0.74);
\draw[thick] (6.5, 0.75) -- (8, 4.25);
\end{pgfonlayer}
  \foreach \x in {0,1,...,2}
     \foreach \y in {0,1,...,2}
        \foreach \z in {0,1,...,1}
     {
\tikzmath{int \c; \c1=0; \c2=\y; \c3=\y+\z; \c4=\x+\y+\z;}
\draw [fill=white,white] (2.75*\x+\y+1.5*\z-0.5, 0.75*\y+3.5*\z-0.2) rectangle (2.75*\x+\y+1.5*\z+0.5, 0.75*\y+3.5*\z+0.2);
\node at (2.75*\x+\y+1.5*\z, 0.75*\y+3.5*\z)  {$\c1\c2\c3\c4$};
}
\end{tikzpicture}
\begin{tikzpicture}[xscale=0.7, yscale=0.7]
\begin{pgfonlayer}{bg}    
\draw[thick] (2,1.5) -- (3.5,5);
\draw[thick] (4.75,1.5) -- (6.25,5);
\draw[thick] (1,0.75) -- (2.5, 4.25);
  \foreach \x in {0,1,...,2}
      \foreach \z in {0,1,...,1}
          \draw[thick] (2.75*\x+1.5*\z, 3.5*\z) -- (2.75*\x+2+1.5*\z, 1.5+3.5*\z);
  \foreach \y in {0,1,...,2}
      \foreach \z in {0,1,...,1}
      {
          \draw[white, fill=white] (\y+1.5*\z-0.1, 0.75*\y+3.5*\z-0.1) rectangle (5.5+\y+1.5*\z+0.1, 0.75*\y+3.5*\z+0.1);
          \draw[thick] (\y+1.5*\z, 0.75*\y+3.5*\z) -- (5.5+\y+1.5*\z, 0.75*\y+3.5*\z);
}           
\draw[thick] (0,0) -- (1.5, 3.5);
\draw[thick] (7.5,1.5) -- (9, 5);
\path[fill=white] (5.4, -0.1) -- (5.6, -0.1) -- (7.1,3.6) -- (6.9,3.6) -- (5.4, -0.1);
\draw[thick] (5.5, 0) -- (7, 3.5);
\draw[thick] (0,0) -- (1.5, 3.5);
\draw[thick] (7.5,1.5) -- (9, 5);
\path[fill=white] (5.4, -0.1) -- (5.6, -0.1) -- (7.1,3.6) -- (6.9,3.6) -- (5.4, -0.1);
\draw[thick] (5.5, 0) -- (7, 3.5);
\path[fill=white] (2.65, -0.1) -- (2.75, -0.1) -- (4.35,3.6) -- (4.15, 3.6) -- (2.65,  -0.1);
\draw[thick] (2.75, 0) -- (4.25, 3.5);
\path[fill=white] (6.4, 0.74) -- (6.6, 0.74) -- (8.1,4.35) -- (7.9,4.35) -- (6.4, 0.74);
\draw[thick] (6.5, 0.75) -- (8, 4.25);
\end{pgfonlayer}
     \foreach \y in {0,1,...,2}
        \foreach \z in {0,1,...,1}
     {
\tikzmath{int \c; \c1=0; \c2=\y; \c3=\y+\z; \c4=1+\y+\z;}
\draw[fill=white,white] (2.75+\y+1.5*\z-0.2, 0.75*\y+3.5*\z-0.2) rectangle (2.75+\y+1.5*\z+0.2, 0.75*\y+3.5*\z+0.2);
\node at (2.75+\y+1.5*\z, 0.75*\y+3.5*\z)  {\scalebox{1.5}{$c$}};
}
     \foreach \y in {1,...,2}
        \foreach \z in {0,...,1}
     {
\tikzmath{int \c; \c1=0; \c2=\y; \c3=\y+\z; \c4=1+\y+\z;}
\draw[fill=white,white] (\y+1.5*\z-0.2, 0.75*\y+3.5*\z-0.2) rectangle (\y+1.5*\z+0.2, 0.75*\y+3.5*\z+0.2);
\node at (\y+1.5*\z, 0.75*\y+3.5*\z)  {\scalebox{1.5}{$c$}};
}
     \foreach \y in {0,...,1}
        \foreach \z in {0,...,1}
     {
\tikzmath{int \c; \c1=0; \c2=\y; \c3=\y+\z; \c4=1+\y+\z;}
\draw[fill=white,white] (5.5+\y+1.5*\z-0.2, 0.75*\y+3.5*\z-0.2) rectangle (5.5+\y+1.5*\z+0.2, 0.75*\y+3.5*\z+0.2);
\node at (5.5+\y+1.5*\z, 0.75*\y+3.5*\z) {\scalebox{1.5}{$c$}};
}
\draw[fill=white,white] (-0.2,-0.2) rectangle (0.2,0.2);
\node at (0,0)  {\scalebox{1.5}{$a$}};
\draw[fill=white,white] (1.3,3.3) rectangle (1.7,3.7);
\node at (1.5, 3.5)  {\scalebox{1.5}{$b$}};
\draw[fill=white,white] (7.5-0.2, 1.5-0.2) rectangle (7.5+0.2, 1.5+0.2);
\node at (7.5, 1.5)  {\scalebox{1.5}{$b$}};
\draw[fill=white,white] (9-0.2, 5-0.2) rectangle (9+0.2, 5+0.2);
\node at (9, 5)  {\scalebox{1.5}{$a$}};
\end{tikzpicture}
}
\caption{On the left, we show $\Pi(0000,0235)/\ZZ\One_4$, from Example~\ref{4Implies1}. On the right, we show the symmetric functions from Example~\ref{2DoesNotImply1}.} \label{NotStrong}
\end{figure}
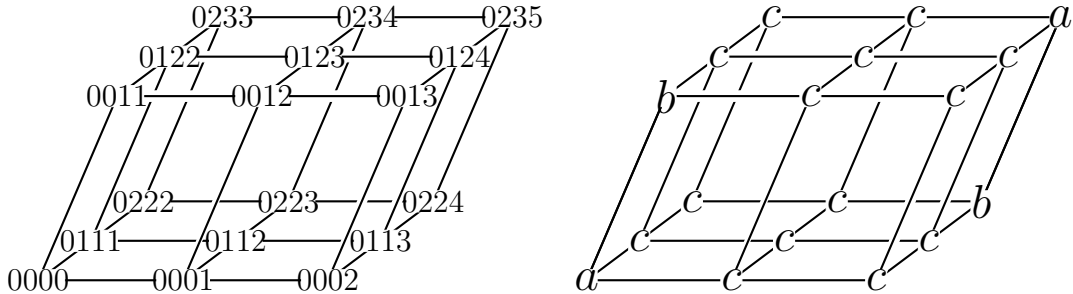

\begin{remark} \label{SchurLlogConcaveRemark}
Let $\Lambda$ denote the ring of symmetric polynomials and write $f \succeq^s g$ to indicate that $f-g$ is Schur nonnegative.
We now discuss the prospects of developing a  theory of Schur $\Llog$-concavity for functions $f: \ZZ^n \to \Lambda$.
For simplicity, assume that $f(x+\One_n) = f(x)$.

There are Schur analogues of all of conditions (1)-(4) in Theorem~\ref{EquivalentLConvexFunctionDefinitions}. Explicitly:
\begin{enumerate}
\item[(1S)]   For any $x$, $y$, $x'$, $y' \in \ZZ^n$ with $x+y=x'+y'$ and $x' \in \Pi(x,y)$, we have \[ f(x) f(y) \preceq^s f(x') f(y') .\]
\item[(2S)] For any $x$, $y \in \ZZ^n$, we have \[ f(x) f(y) \preceq^s f(x \inmin y) f(x \inmax y) .\]
\item[(3S)] For any $x$, $y \in \ZZ^n$, we have \[ f(x) f(y) \preceq^s f(\floorAve{x}{y}) f(\ceilAve{x}{y}) . \]
\item[(4S)] The set $\Supp(f):= \{ x\in \ZZ^n : f(x) \neq 0 \}$ is $L$-convex and, for any $x \in \ZZ^n$ and any subsets $I \subseteq J$ of $[n]$, we have \[ f(x)  f(x+e_I+e_J) \preceq^s f(x+e_I)  f(x+e_J) . \]
\end{enumerate}
The proofs of $(1S) \! \implies \! (2S) \! \implies \! (3S) \! \implies \! (4S)$ proceed exactly as in Theorem~\ref{EquivalentLConvexFunctionDefinitions}.

However, the implication from (4S) to (1S) is no longer valid. The issue is that, in our proof of the implication $(4) \implies (1)$, we repeatedly use the fact that, if $f \in \RR$ and $f+g \geq f+h$, then $g \geq h$. For example, in Example~\ref{4Implies1}, we must cancel six terms from each side of an inequality. If we were to try to mimic this proof for Schur positivity, we would need it to be true that, if $f\neq 0$ is Schur nonnegative and $fg \succeq^s fh$, then $g \succeq^s h$. This is not valid.

Lam, Postnikov and Pylyavskyy~\cite{LPP} have already proved that \[ s_{\floorAve{\lambda}{\mu}} s_{\ceilAve{\lambda}{\mu}} \succeq^s s_{\lambda} s_{\mu} \ \text{and}\ s_{\lambda \inmin \mu} s_{\lambda \inmax \mu}  \succeq^s s_{\lambda} s_{\mu}. \] 
(These inequalities were conjectured in~\cite{Oko1} and~\cite{FFLP}.)
Thus, ignoring the condition about $\One_n$, they prove that $\lambda \mapsto s_{\lambda}$ is Schur-$\Llog$-concave in the sense of conditions (2S) or (3S).
However, these are not enough to imply condition~(1S) and thus prove the LPP conjecture.

If we take $n=2$ and work with functions which are invariant under translation by $\One_2$, then this is equivalent to working with functions on $\ZZ^2/\ZZ \One_2 \cong \ZZ$. In this guise of functions $\ZZ \to \Lambda$,  Matherne, Miyata, Proudfoot and Ramos~\cite{Equivariant} define a function satisfying~(1S) to be ``strongly equivariantly $\log$-concave" and define a function satisfying~(4S) to be ``weakly equivariantly $\log$-concave".
The author believes that the strong condition, (1S), is the correct condition to consider.
\end{remark}

\begin{example} \label{2DoesNotImply1}
We give an example of a function $f:\ZZ^4 \to \Lambda$ which obeys condition~(2S)  but not~(1S).
We give the values of $f$ on $\{ 0 \} \times \ZZ^3$, and extend by $\One_4$ invariance.
See Figure~\ref{NotStrong}.
\[
\begin{array}{r@{}r@{}l@{}c@{}l}
f_{0000} = f_{0235} &=& a &:=& s_{21} \\
f_{0011} = f_{0224} &=& b &:=& 2 s_3 + 2 s_{111} \\
 \begin{matrix} f_{0001}=f_{0002}=f_{0112}=f_{0113}=f_{0222}=f_{0223}=f_{0224} \\ =f_{0011}=f_{0012}=f_{0013}=f_{0122}=f_{0123}=f_{0124}=f_{0233} \end{matrix} 
 &\left. \vphantom{\begin{matrix} f_{0000} \\ f_{0000} \end{matrix}} \right\} \! =& c &:=&  2 s_3 + 4 s_{21} + 2 s_{111}  \\
f_{0xyz} &=& 0 && \ \text{otherwise} 
\end{array}
\]
Then $f_{0000} f_{0235} = s_{21}^2$ contains an $s_{321}$ term, whereas $f_{00011} f_{0224} = ( 2 s_3 + 2 s_{111})^2$ does not, so $f_{0000} f_{0235} \not\preceq^s f_{0011} f_{0224}$. However,  $f_{\lambda} f_{\mu} \preceq^s f_{\lambda \inmin \mu} f_{\lambda \inmax \mu}$ for all $\lambda$, $\mu \in \ZZ^4$.
\end{example}

\section{Marginalizations of $\Llog$-concave functions are $\Llog$-concave} \label{LConcaveProof}

The point of this section is to prove Theorem~\ref{PointCountLConcave}.

\begin{theorem} \label{PointCountLConcave}
Let $K \subset \ZZ^{M+N}$ be an $L$-convex set. For $x \in \ZZ^M$, let \[E_K(x) = \#\{ y \in \ZZ^N : (x,y) \in K \},\]and assume that $E_K(x) < \infty$ for all $x \in \ZZ^M$. Then $E_K$ is an $\Llog$-concave function on $\ZZ^M$.
\end{theorem}

Here, if $x = (x_1, \ldots, x_n)$ and $y = (y_1, \ldots, y_m)$, then the notation $(x,y)$ is shorthand for $(x_1, \ldots, x_n, y_1, \ldots, y_m)$.

Our starting point is to turn $L$-convexity of sets into $L$-convexity of functions:
\begin{lemma} \label{CharacteristicLConvex}
Let $K \subset \ZZ^N$ be an $L$-convex set, and let $\chi_K : \ZZ^N \to \{ 0, 1 \}$ be the characteristic function of $L$. Then $\chi_K$ is $\Llog$-concave.
\end{lemma}

\begin{proof}
Invariance under translation by $\One_N$ is clear.

Let $x$, $y$, $x'$ and $y'$ be as in the definition of $L$-convexity (Definition~\ref{LConvexFunctionDefn}); we must show that $f(x) f(y) \leq f(x') f(y')$. 
If $x$ or $y \not\in K$, then the left hand side is $0$ and we are done. If $x$ and $y \in K$, then we have $x'$ and $y' \in K$ because $K$ is $L$-convex, so our desired inequality is $1 \cdot 1 \leq 1 \cdot 1$ and we are done.
\end{proof}

Our proof is based on the Ahlswede–Daykin inequality. See~\cite{AD} for the original proof and see~\cite[Chapter 6]{AlonSpencer} for modern applications.

Let $X$ be a distributive lattice, meaning a set with two commutative, associative, idempotent binary options $\inmin$ and $\inmax$, each of which distributes over the other, and which obey the absorption identities: $x \inmin (x \inmax y) = x$ and $x \inmax (x \inmin y) = x$. 
In our setting, $X$ will be $\ZZ^N$, and these will be our standard operations $\inmin$ and $\inmax$.

 Let $U$ and $V$ be subsets of $X$ and define $U \inmin V := \{ u \inmin v : u \in U,\ v \in V \}$ and $U \inmax V := \{ u \inmax v : u \in U,\ v \in V \}$. 
Let $f_1$, $f_2$, $f_3$ and $f_4$ be $\RR_{\geq 0}$-valued functions on $U$, $V$, $U \inmin V$ and $U \inmin V$ respectively. 

\begin{AD}
Suppose that, for all $u \in U$ and $v \in V$, we have $f_1(u) f_2(v) \leq f_3(u \inmin v) f_4(u \inmax v)$. Then 
\[ \left( \sum_{u \in U} f_1(u) \right)  \left( \sum_{v \in V} f_1(v) \right) \leq \left( \sum_{t \in U \inmin V} f_3(t) \right)  \left( \sum_{w \in U \inmax V} f_4(w) \right) . \]
\end{AD}

We now begin our proof.

\begin{Lemma} \label{PushforwardLConcave}
Let $h : \ZZ^{M+N} \to \RR_{\geq 0}$ be an $\Llog$-concave function with $h(x+\One_{M+N}) = R h(x)$. Define $\bar{h} : \ZZ^{M} \to \RR_{\geq 0}$ by
\[ \bar{h}(x_1, x_2, \ldots, x_M) = \sum_{(y_1, y_2, \ldots, y_N) \in \ZZ^N} h(x_1, x_2, \ldots, x_M,  y_1, y_2, \ldots, y_N) \]
and assume that the sum defining $\bar{h}$ is always finite. Then $\bar{h}$ is an $\Llog$-concave function $\ZZ^{M} \to \RR_{\geq 0}$, obeying $\bar{h}(x+\One_{M}) = R\bar{h}(x)$. 
\end{Lemma}


\begin{proof}
We first check the criterion that $\bar{h}(x+\One_{M}) = R\bar{h}(x)$. Indeed,
\[ \bar{h}(x+\One_{M}) = \sum_{y \in \ZZ^N} h(x+\One_{M}, y) =  \sum_{y \in \ZZ^N}  (R h(y, y-\One_N)) = R \bar{h}(x). \]

We now check condition~(2) for $\Llog$-concavity. We must show, for $x$, $y \in \ZZ^{M}$, that $\bar{h}(x) \bar{h}(y) \leq \bar{h}(x \inmin y) \bar{h}(x \inmax y)$.
In other words, we must show that
\[ 
\left( \sum_{p \in \ZZ^N} h(x, p) \right) \left( \sum_{q \in \ZZ^N} h(y,q) \right) \leq \left( \sum_{r \in \ZZ^N} h(x \inmin y, r) \right) \left( \sum_{s \in \ZZ^N} h(x \inmax y, s)  \right).
\]
This is precisely the Ahlswede–Daykin inequality, with $f_1=f_2=f_3=f_4$, and $U = \{ (x,p) : p \in \ZZ^N \}$ and $V = \{ (y,q) : q \in \ZZ^N \}$. 
\end{proof}

Theorem~\ref{PointCountLConcave} now follows by applying Lemma~\ref{PushforwardLConcave} to the function $\chi_K$, which is $\Llog$-concave by Lemma~\ref{CharacteristicLConvex}.

\begin{remark} \label{ThankRutherford}
The author asked for a proof of Lemma~\ref{PushforwardLConcave} on Mathoverflow~\cite{MO504682} and received a response from user ``Ambrose Rutherford" that was likely in violation of Mathoverflow's policies on AI generated content; as of current writing, Ambrose Rutherford's reply is deleted.
The author nonetheless found this response extremely useful, as it told him of the Ahlswede–Daykin inequality, which he suspects it would have taken him weeks if not months to find on his own.
\end{remark}
%
%

\section{Proofs of Theorems~\ref{SkepExtLConcave}, \ref{BetterLPP} and the LPP conjecture} \label{Finishing}

We are now ready to prove our main theorems. 

\begin{Theorem} \label{SkepExtLConcave}
Fix $g^+ : \Delta^+_n \to \ZZ$. Then $\SkepExt(g^+, \rho)$ is $\Llog$-concave as a function of $\rho$.
\end{Theorem}

\begin{proof}
Without loss of generality, we may add the same constant to every element of $g^+$ in order to assume that $g^+_{n0} =0$. 

Fixing $g^+$, the condition that $(g^+, g^-)$ is a skep is a collection of inequalities of the form $g^-_{ij} - g^-_{k \ell} \leq c_{ij,k\ell}$. 
In other words, the set of $g^-$ such that $(g^+, g^-)$ is a skep is $L$-convex. 
We will apply Theorem~\ref{PointCountLConcave} to this $L$-convex set, projecting onto the coordinates $(g^-_{(n-1)1}, g^-_{(n-3)3}, \ldots, g^-_{1(n-1)})$.
So the number of $g^-$ such that $(g^+, g^-)$ is a skep and $(g^-_{(n-1)1}, g^-_{(n-3)3}, \ldots, g^-_{1(n-1)}) = \tau$  is $\Llog$-concave as a function of $\tau$.

Set $\sigma := (g^+_{n0}, g^+_{(n-2)2}, \ldots, g^+_{2(n-2)})$.
By definition,
\[ \partial_1(g^+, g^-) =  (g^-_{(n-1)1}, g^-_{(n-3)3}, \ldots, g^-_{1(n-1)}) - (g^+_{n0}, g^+_{(n-2)2}, \ldots, g^+_{2(n-2)}) \]
so $\partial_1(g^+, g^-)  = \rho$ if and only if $(g^-_{(n-1)1}, g^-_{(n-3)3}, \ldots, g^-_{1(n-1)}) = \rho+\sigma$. 
So $\SkepExt(g^+, \rho)$  is $\Llog$-concave as a function of $\rho+\sigma$. 
By Corollary~\ref{Translate}, $\SkepExt(g^+, \rho)$ is also $\Llog$-concave as a function of $\rho$.
\end{proof}

\begin{Theorem} \label{BetterLPP}
Let $\lambda$, $\mu$ and $\nu$ be partitions. Set $\pi = \lambda+\mu$.
Fix $g^+ : \Delta^+_n \to \ZZ$ with $\DB(g^+) = \nu$ and $\OB(g^+) = \pi$.
Let $\lambda'+\mu' = \pi$ with $\lambda' \in \Pi(\lambda, \mu)$. 
Then
\begin{equation}
\SkepExt(g^+, \lambda) \leq \SkepExt(g^+, \lambda').  \label{SkepExtInequality}
\end{equation}
\end{Theorem}

\begin{proof}
Since $\lambda' \in \Pi(\lambda, \mu)$, we have $\mu' \in \Pi(\lambda, \mu)$. So, by Theorem~\ref{SkepExtLConcave}, we have
\begin{equation}
 \SkepExt(g^+, \lambda) \SkepExt(g^+, \mu)  \leq \SkepExt(g^+, \lambda') \SkepExt(g^+, \mu') . \label{SquareOfBetterLPP}
 \end{equation}
By Lemma~\ref{SkepCommutative}, $\SkepExt(g^+, \lambda) =  \SkepExt(g^+, \mu)$ and  $\SkepExt(g^+, \lambda') =  \SkepExt(g^+, \mu')$, so we can rewrite Equation~\eqref{SquareOfBetterLPP} as
\[  \SkepExt(g^+, \lambda)^2  \leq \SkepExt(g^+, \lambda')^2 .\]
Taking square roots, we have the claim.
\end{proof}

\begin{LPPMain}
Let $\lambda$, $\mu$, $\lambda'$, $\mu'$ and $\nu$ be partitions in $\ZZ^n$. Suppose that $\lambda'+\mu' = \lambda+\mu$ and $\lambda' \in \Pi(\lambda, \mu)$.   Then $c_{\lambda \mu}^{\nu} \leq c_{\lambda' \mu'}^{\nu}$. 
\end{LPPMain}

\begin{proof}[Proof of the Main Theorem]
From Corollary~\ref{SkepLRAsSum}, 
\[ c_{\lambda \mu}^{\nu} = \sum \SkepExt(g^+, \lambda)  \ \text{and} \ 
 c_{\lambda' \mu'}^{\nu} = \sum \SkepExt(g^+, \lambda')  \]
where both sums are over the same set: functions $g^+ : \Delta_n^+ \to \ZZ$ obeying $\DB(g^+) = \nu$, $\OB(g^+) = \lambda+\mu$ and $g_{n0}=0$.
From Theorem~\ref{BetterLPP}, the $g^+$ summand in the first term is always $\leq$ the $g^+$ term in the second sum. So $c_{\lambda \mu}^{\nu} \leq c_{\lambda' \mu'}^{\nu'}$
as desired.
\end{proof}

\section{The inequalities necessary to prove the LPP conjecture} \label{MinimalInequalities}
It is natural to approach the LPP conjecture by finding a list of quadruples $(\lambda, \mu, \lambda', \mu')$ such that, if we prove the LPP conjecture for these quadruples, it follows for all others. 
The point of this section is to give such a list. 
Our final proof doesn't use this list, but it was valuable in finding the proof to know which cases to check, so it seems worth recording the list here. 
This list isn't quite minimal (see Remark~\ref{NotMinimal}), but removing the nonminimal elements would make the statement of Theorem~\ref{Covers} significantly longer with little benefit.

\begin{remark} \label{ThankSawin}
Key parts of this result were found by Will Sawin, in  response to the author's question~\cite{MO422985} on MathOverflow. 
More specifically, the author found Lemma~\ref{Contracting}, reformulating the problem in terms of contracting maps, and then Sawin found the minimal contracting maps.
The author appreciates Sawin giving permission for his answer to appear here. 
\end{remark}

\begin{definition}
Let $S_{\pi}$ be the set of ordered pairs of vectors $(\lambda, \mu)$ such that $\lambda+\mu = \pi$. We are explicit that we don't assume $\lambda$ and $\mu$ to be dominant.
 For $(\lambda, \mu)$ and $(\lambda', \mu') \in S_{\pi}$, define $(\lambda', \mu') \sqsubseteq (\lambda, \mu)$ if $\Pi(\lambda', \mu') \subseteq \Pi(\lambda, \mu)$.

The relation $\sqsubseteq$ is a preorder on $S_{\pi}$, meaning that it is reflexive and transitive.
Recall that, whenever we have a preorder $R$ on a set $X$, there is an associated equivalence relation where $x$ is equivalent to $y$ if $x R y R x$.
We write $\sim$ for the equivalence relation associated to $\sqsubseteq$, and we'll write $(\lambda', \mu') \sqsubset (\lambda, \mu)$ if $(\lambda', \mu') \sqsubseteq (\lambda, \mu)$ and  $(\lambda, \mu) \not\sim (\lambda', \mu')$.
\end{definition}

%
%

Thus, the LPP conjecture states that, if $(\lambda, \mu) \sqsupseteq (\lambda', \mu')$, then $s_{\lambda'} s_{\mu'} - s_{\lambda} s_{\mu}$ is Schur nonnegative. 
The point of this section is to give a list of relations $(\lambda, \mu) \sqsupset (\lambda', \mu')$ such that, whenever we have $(\lambda, \mu) \sqsupset (\lambda'', \mu'')$, we can find a relation in our list such that $(\lambda, \mu) \sqsupset (\lambda', \mu') \sqsupseteq (\lambda'', \mu'')$.
Since $\Pi(\lambda, \mu)/\ZZ \One_n$ is finite, if one proves the LPP conjecture for all $(\lambda, \mu, \lambda', \mu')$ in our list, then one has proved the full LPP conjecture.

\begin{definition}
For integers $b<c$, and other integers $x$ and $y$, we define
\[ x \zig_{bc} y = 
\begin{cases} 
x+b & y \leq x+b \\
y & x+b \leq y \leq x+c \\
x+c &  x+c \leq y \\
\end{cases} \qquad \qquad
x \zag_{bc} y = 
\begin{cases} 
y-b &  y \leq x+b \\
x & x+b \leq y \leq x+c \\
y-c & x+c \leq y.
\end{cases}\]
We pronounce $\zig$ as ``zig" and $\zag$ as ``zag".
For $\lambda$ and $\mu \in \ZZ^n$, we define $\lambda \zig_{bc} \mu$ and $\lambda \zag_{bc} \mu$ by applying the zig and zag operators in each coordinate. 
\end{definition}

\begin{remark}
We have $x \zig_{0 \infty} y = x \inmax y$ and $x \zag_{0 \infty} y = x \inmin y$, whereas $x \zig_{(-\infty)0} y = x \inmin y$ and $x \zig_{(-\infty)0} y = x \inmax y$. Thus, zig and zag generalize min and max.
\end{remark}

We observe that $x \zig_{bc} y + x \zag_{bc} y = x+y$ so, if $(\lambda,\mu) \in S_{\pi}$, then $(\lambda \zig_{bc} \mu,\ \lambda \zag_{bc} \mu) \in S_{\pi}$ as well. 
It is also easy to check that $(\lambda, \mu) \sqsubseteq (\lambda \zig_{bc} \mu,\ \lambda \zag_{bc} \mu)$, and that $(\lambda, \mu) \sqsubset (\lambda \zig_{bc} \mu,\ \lambda \zag_{bc} \mu)$ as long as either $b$ or $c$ lies in $[\min(\delta_i), \max(\delta_i)]$.

\begin{theorem} \label{Covers}
For any $(\lambda, \mu) \sqsupseteq (\lambda'', \mu'')$ in $S_{\pi}$, we can find $b < c$ such that \[ (\lambda, \mu) \sqsupset (\lambda \zig_{bc} \mu, \lambda \zag_{bc} \mu) \sqsupseteq (\lambda'', \mu'').\] Moreover, we can take $(b,c)$ to be of the form either $(\delta_p, \delta_{p}+1)$ or $(\delta_p, \delta_q)$.
\end{theorem}

%

\begin{remark} \label{NotMinimal}
In general, the relations in Theorem~\ref{Covers} are not a minimal list.
Let $(\lambda, \mu) = (1111, 1234)$ and let $(\lambda', \mu') = (1111 \zig_{12} 1234,\ 1111 \zag_{12} 1234) = (2233, 0112)$.
Then \[(1111, 1234) \sqsupset (1112, 1233) \sqsupset (2233, 0112)\] so  the relation  $(\lambda, \mu) \sqsupset (\lambda', \mu')$ is redundant.
\end{remark}

\begin{remark} \label{ZigZagGeometry}
In terms of the parallelepiped $\Pi(\lambda, \mu)/\ZZ \One_n$, the points $\lambda \zig_{\delta_p (\delta_p+1)} \mu$ are neighbors of $\lambda$ along the edges of the parallelepiped, so they are ``nearby" $\lambda$.
The points $\lambda \zig_{\delta_p \delta_q} \mu$, on the other hand, are vertices of the parallelepiped and can be ``far from" $\lambda$ and $\mu$. 
\end{remark}

\begin{eg} \label{Example0358}
Let $\lambda$ be some partition with $4$ rows and let $\mu = \lambda + (0,3,5,8)$. Then $\Pi(\lambda, \mu) / \ZZ \One_4$ has $48$ elements, which are depicted in Figure~\ref{Example0358Figure}. For legibility, we have taken $\lambda = (0,0,0,0)$ in the figure, but it would be more interesting from the perspective of the LPP conjecture if $\lambda$ were large, since the LPP conjecture is trivial when $\lambda = 0$.

The elements in rectangular boxes correspond to covers of the first type, and the elements in oval boxes correspond to covers of the second type.
In this example, all of these are genuine covers, so any proof of the LPP conjecture must handle these cases directly, and not as a consequence of other cases. In particular, $\lambda+(0,0,2,2)$ and $\lambda+(0,3,3,6)$ are neither neighbors of $\lambda$ or $\mu$, nor are they of the form $(\lambda+a \One_4) \inmin (\mu+b \One_4)$ or $(\lambda+a \One_4) \inmax (\mu+b \One_4)$.
\end{eg}

\usetikzlibrary{math}
\usetikzlibrary{shapes}

\begin{figure}
\centerline{
\begin{tikzpicture}[scale=0.75]
\begin{pgfonlayer}{bg}    
\draw[thick] (0,0) -- (3,7) ;
\draw[thick, dashed] (3,2.25) -- (6,9.25) ;
\foreach \x in {0,1,...,3}
    \foreach \z in {0,1,...,2}
          \draw[thick] (2.75*\x+1.5*\z, 3.5*\z) -- (2.75*\x+3+1.5*\z, 2.25+3.5*\z);
  \foreach \y in {0,1,...,3}
      \foreach \z in {0,1,...,2}
      {
          \draw[white, fill=white] (\y+1.5*\z-0.1, 0.75*\y+3.5*\z-0.1) rectangle (8.25+\y+1.5*\z+0.1, 0.75*\y+3.5*\z+0.1);
          \draw[thick] (\y+1.5*\z, 0.75*\y+3.5*\z) -- (8.25+\y+1.5*\z, 0.75*\y+3.5*\z);
}           
\draw[thick] (0,0) -- (3, 7);
\draw[thick] (11.25,2.25) -- (14.25, 9.25);
\path[fill=white] (8.15, -0.1) -- (8.35, -0.1) -- (11.35, 7.1) -- (11.15, 7.1) -- (8.15, -0.1);
\draw[thick] (8.25, 0) -- (11.25, 7);
\end{pgfonlayer}

\node[draw, fill=white] at (2.75, 0)  {$\phantom{0001}$};  
\node[draw, fill=white] at (1.5, 3.5)  {$\phantom{0011}$}; 
\node[draw, fill=white] at (1, 0.75)  {$\phantom{0111}$};  
\node[draw, fill=white] at (11.5, 9.25)  {$\phantom{0357}$};  
\node[draw, fill=white] at (12.75, 5.75)  {$\phantom{0347}$};
\node[draw, fill=white] at (13.25, 8.5)  {$\phantom{0247}$}; 

\node[draw, shape=ellipse, fill=white] at (8.25, 0)  {$\phantom{0003}$}; 
\node[draw, shape=ellipse, fill=white] at (3, 7)  {$\phantom{0022}$}; 
\node[draw, shape=ellipse, fill=white] at (3, 2.25)  {$\phantom{0333}$}; 
\node[draw, shape=ellipse, fill=white] at (6, 9.25)  {$\phantom{0355}$}; 
\node[draw, shape=ellipse, fill=white] at (11.25, 2.25)  {$\phantom{0336}$}; 
\node[draw, shape=ellipse, fill=white] at (11.25, 7)  {$\phantom{0025}$};

  \foreach \x in {0,1,...,3}
     \foreach \y in {0,1,...,3}
        \foreach \z in {0,1,...,2}
     {
\tikzmath{int \c; \c1=0; \c2=\y; \c3=\y+\z; \c4=\x+\y+\z;}
\draw [fill=white,white] (2.75*\x+\y+1.5*\z-0.5, 0.75*\y+3.5*\z-0.2) rectangle (2.75*\x+\y+1.5*\z+0.5, 0.75*\y+3.5*\z+0.2);
\node at (2.75*\x+\y+1.5*\z, 0.75*\y+3.5*\z)  {$\c1\c2\c3\c4$};
}

\end{tikzpicture}}
\caption{The parallelepiped $\Pi(0000,\ 0358)/\ZZ \One_4$, discussed in Example~\ref{Example0358}.} \label{Example0358Figure}
\end{figure}

We now introduce a helpful change of coordinates.
\begin{definition}
Let $C_{\pi}$ be the set of $\delta \in \ZZ^n$ such that $\delta \equiv \pi \bmod 2$. 
Put a preorder $\preceq$ on $C_{\pi}$ by $\delta'  \preceq \delta$ if for all $1 \leq i<j \leq n$, we have $|\delta'_i - \delta'_j| \leq |\delta_i - \delta_j|$. 

We write $\equiv$ for the equivalence relation corresponding to the preorder $\preceq$. Concretely, $\delta \equiv \delta'$ if $\delta'$ is of the form $\pm \delta + c \One_n$ for some $c \in \ZZ$. We write $\delta' \prec \delta$ if $\delta' \preceq \delta$ and $\delta \not\equiv \delta'$.
\end{definition}

\begin{lemma}
The map $(\lambda, \mu) \mapsto \mu - \lambda$ is an isomorphism of preorders between $S_{\pi}$ and $C_{\pi}$.
The inverse map is $\delta \mapsto (\tfrac{\pi-\delta}{2},\ \tfrac{\pi+\delta}{2})$.
\end{lemma}


\begin{proof}
The two maps are clearly mutually inverse bijections between $S_{\pi}$ and $C_{\pi}$. We need to check that the maps respect the preorders. 
Let $(\lambda, \mu)$, $(\lambda', \mu') \in S_{\pi}$ and put $\delta = \mu-\lambda$, $\delta' = \mu' - \lambda'$. By definition, $(\lambda, \mu) \sqsupseteq (\lambda', \mu')$ if, for all $1 \leq i < j \leq n$, we have
\begin{equation} \min(\lambda_i-\lambda_j, \mu_i-\mu_j) \leq \lambda'_i-\lambda'_j \leq \max(\lambda_i-\lambda_j, \mu_i-\mu_j) . \label{blah} \end{equation}
We convert each part of~\eqref{blah} into formulas in $\pi$, $\delta$ and $\delta'$:
\[ \begin{array}{lcl}
 \min(\lambda_i - \lambda_j, \mu_i - \mu_j) 
 &=& \tfrac{1}{2} \left( \pi_i - \pi_j - |\delta_i - \delta_j| \right) \\[0.2 cm]
 \lambda'_i -  \lambda'_j
  &=& \tfrac{1}{2} (\pi_i -\pi_j -\delta'_i + \delta'_j) \\[0.2 cm]
 \max(\lambda_i - \lambda_j, \mu_i - \mu_j)  
  &=& \tfrac{1}{2} \left( \pi_i - \pi_j + |\delta_i - \delta_j| \right) .
 \end{array} \]
Cancelling common terms,~\eqref{blah} holds if and only if
\[ -|\delta_i - \delta_j| \leq -\delta'_i + \delta'_j \leq |\delta_i-\delta_j| \]
or, in other words, $|\delta'_i - \delta'_j| \leq |\delta_i-\delta_j|$.
 \end{proof}

We now give another, more geometric, description of the preorder $\succeq$.
\begin{definition}
Define  $\phi: \RR \to \RR$ to be a \newword{contraction} if $|\phi(x)-\phi(y)| \leq |x-y|$.
\end{definition}

\begin{lemma} \label{Contracting}
We have $\delta \succeq \delta'$ if and only if we can find a contraction $\phi$ with $\phi(\delta_i) = \delta'_i$.
\end{lemma}

\begin{proof}
If we can find such a $\phi$, then we have $|\delta'_i - \delta'_j| = |\phi(\delta_i) - \phi(\delta_j)| \leq |\delta_i - \delta_j|$ so $\delta' \preceq \delta$. 

Conversely, suppose that $\delta' \preceq \delta$. First, note that, if $\delta_i = \delta_j$ then $\delta'_i = \delta'_j$, so we can find a function $\phi: \RR \to \RR$ mapping the $\delta_i$'s to the $\delta'_i$'s. If we take this function to be piecewise linear with corners at the $\delta_i$, then the slopes of $\phi$ are $\tfrac{\delta'_i - \delta'_j}{\delta_i - \delta_j}$ for $\delta_i$ and $\delta_j$ consecutive.
The assumption that $\delta' \preceq \delta$ implies that these slopes are in $[-1,1]$, so $\phi$ is a contraction.
\end{proof}

The zig/zag operators correspond to specific contractions.
Define
\[ \phi_{bc}(z) = \begin{cases} z-2b & z \leq b \\ -z & b \leq z \leq c \\ z-2c & c \leq z \end{cases} .\]
We remark that the graph of $\phi$ resembles $\zag$. An immediate computation shows:

\begin{lemma}
Let $(\lambda, \mu) \in S_{\pi}$ and put $\delta = \mu-\lambda$. Let $(\lambda', \mu') = (\lambda \zig_{bc} \mu,\ \lambda \zag_{bc} \mu)$ and put $\delta' = \mu' - \lambda'$. Then $\delta'_i = \phi_{bc}(\delta_i)$.
\end{lemma}


\begin{proof}[Proof of Theorem~\ref{Covers}]
Let $\delta  \succ \delta''$ in $C_{\pi}$. 
We must find $(b,c)$ of the form either $(\delta_p, \delta_p +1)$ or $(\delta_p, \delta_q)$, such that $\delta \succ \phi_{bc}(\delta) \succeq \delta''$.
Rename our coordinates so that $\delta_1 \leq \delta_2 \leq \cdots \leq \delta_n$.

\textbf{Case 1:} There is some index $p$ such that $|\delta''_{p+1} - \delta''_p| < \delta_{p+1} - \delta_p$. Note that, since $\delta''_i \equiv \delta_i \bmod 2$, we have  $|\delta''_{p+1} - \delta''_p| \leq \delta_{p+1} - \delta_p-2$.

 Abbreviate $\delta_p$ to $d$ and put $\delta':= \phi_{d(d+1)}(\delta)$. We will show that $\delta' \succeq \delta''$. We have $\delta'_i = \delta_i - 2d$ for $i \leq p$ and $\delta'_i = \delta_i - 2d-2$ for $i \geq p+1$, so $\delta'_{i+1} - \delta'_i = \delta_{i+1} - \delta_i$ for $i \neq p$ and $\delta'_{p+1} - \delta'_p = \delta_{p+1} - \delta_p-2$. Thus, $|\delta''_{i+1}-\delta''_i| \leq \delta'_{i+1} - \delta'_i$ for all $i$. Since $\delta'_1 \leq \delta'_2 \leq \cdots \leq \delta'_n$, this forces $|\delta''_{j}-\delta''_i| \leq \delta'_{j} - \delta'_i$ for all $1 \leq i<j \leq n$, and hence $\delta'' \preceq \delta'$, as claimed. 

\textbf{Case 2:} For all $1 \leq i < n$, we have  $|\delta''_i - \delta''_{i+1}| = |\delta_i - \delta_{i+1}|$. 
Let $\psi : [\delta_1, \delta_n] \to \RR$ be the unique contraction with $\psi(\delta_i) = \delta''_i$, so all slopes of $\psi$ are $\pm 1$. 
Define an interval $[\delta_i, \delta_j]$ to be a \newword{segment} of $\psi$ if $\psi$ is affine linear on $[\delta_i, \delta_j]$, but not on any larger interval containing $[\delta_i, \delta_j]$. 
If $\psi$ has only one segment, then $\delta \equiv \delta''$, contrary to our hypothesis that $\delta \succ \delta''$.
So we may assume that $\psi$ has at least two segments.

Let $[\delta_p, \delta_q]$ be a segment of minimal length. We abbreviate $\delta_p$ and $\delta_q$ to $b$ and $c$. 
Let $[a,b]$ be the segment to the left of $[b,c]$; if $b=\delta_1$, we put $a = -\infty$.
Let $[c,d]$ be the segment to the right of $[b,c]$; if $c=\delta_n$, we put $d = \infty$.
Put $\delta'' = \phi_{b c}(\delta)$. We must show that $\delta'' \succeq \delta'$.

Define a contraction $\theta$ by
\begingroup 
\setlength{\belowdisplayskip}{2pt} \setlength{\belowdisplayshortskip}{2pt}
\setlength{\abovedisplayskip}{2pt} \setlength{\abovedisplayshortskip}{2pt}
\begin{equation}
\theta(x) := \begin{cases} \psi(x+2b) & x \leq -b \\ \psi(x+2c) & x \geq -c \end{cases}. \label{psiDefn}
 \end{equation}
 \endgroup
Since $b<c$, we have $-b>-c$, and thus~\eqref{psiDefn} assigns a value to $\theta(x)$ for all $x$. We need to check that, if $x \in [-c, -b]$, then the two formulas in~\eqref{psiDefn} are equal.
We do this now.

Without loss of generality, assume that $\psi$ has slopes $1$, $-1$ and $1$ on $[a,b]$, $[b,c]$ and $[c,d]$ respectively.
If $x \in [-c, -b]$, then $x+2b \in [2b-c, b] \subseteq [a,b]$, where we have used that $[b,c]$ is the shortest segment to deduce that $b-a \geq c-b$ and hence $2b-c \geq a$. 
Since $\psi$ has slope $1$ on $[a,b]$, we have $\psi(x+2b) = x+b+\psi(b)$. Similarly, $\psi(x+2c) =x+c+\psi(c)$. Since $\psi$ has slope $-1$ on $[b,c]$, we have $b+\psi(b) = c+\psi(c)$, and the two formulas are compatible on the region of overlap. 
Since $\psi$ is a contraction on $(\infty, b]$ and is a contraction on $[c, \infty)$, the function $\theta$ is a contraction on $(\infty, -b]$ and $[-c, \infty)$, and hence $\theta$ is a contraction.

We now verify that $\theta(\delta'_i) = \delta''_i$, thus confirming that $\delta' \succeq \delta''$ by Lemma~\ref{Contracting}.
If $i \leq p$, then $\delta'_i = \phi_{bc}(\delta_i) = \delta_i - 2b \leq \delta_p - 2b = -b$ so $\theta(\delta'_i) = \psi(\delta'_i+2b) = \psi(\delta_i) = \delta''_i$. 
If $i \geq p+1$ then $\delta'_i = \phi_{bc}(\delta_i) = \delta_i - 2c \geq \delta_q -2c = -c$ so $\theta(\delta'_i) = \psi(\delta'_i+2c) = \psi(\delta_i) = \delta''_i$. 
We have proven that  $\theta(\delta'_i) = \delta''_i$ and thus  $\delta' \succeq \delta''$. 
\end{proof}

\vspace{-0.25 in}

\section*{Acknowledgements} 
The author has discussed a different approach to the LPP conjecture with Jake Levinson many times, and is still interested in pursuing it.

Alex Postnikov's talk at Open Problems in Algebraic Combinatorics changed the LPP conjecture from something that the author was somewhat aware of to a concrete goal; the author thanks Alex for the talk and OPAC for its excellent atmosphere for discussions.

The invention of skeps was spurred by a conversation with Son Nguyen where Nguyen noted that, in most combinatorial models for $c_{\lambda \mu}^{\nu}$, the term $\lambda_i$ is ``far from" $\mu_i$, making it hard to use the hypothesis $\lambda+\mu = \lambda' + \mu'$. Nguyen and Pylyavskyy's ``shuffle tableaux"~\cite{NguyenPylavskyy, NguyenNguyenWoodruff} are another attempt to solve this problem. The author suspects that there are deeper relations between skeps and shuffle tableaux.


The author thanks Mathoverflow users  Ambrose Rutherford and Will Sawin, as described in Remarks~\ref{ThankRutherford} and~\ref{ThankSawin}.
The author thanks Tuong Le for pointing out that the issue in Corollary~1 was addressed incorrectly in an earlier version.
ChatGPT told the author that Condition~(0) in Theorem~\ref{EquivalentLConvexSetDefinitions} was called ``L-convexity", and provided a mostly correct proof of Theorem~\ref{PointCountLConcave} for $M=2$, but it also told the author that  Theorem~\ref{PointCountLConcave} was false for $M \geq 3$.
The author thanks Allen Knutson, Kazuo Murota and Pavlo Pylavskyy for helpful comments on this manuscript.

The author was supported in part by NSF grant DMS-2246570.

\raggedright

\end{document}